\documentclass[10pt,journal,twocolumn,oneside,letterpaper,final]{IEEEtran}

\usepackage{amsmath,amssymb,amsfonts,amsthm}
\usepackage{dsfont,mathpazo,mathrsfs}
\usepackage{esvect}

\usepackage{graphicx,tikz} 
\usepackage[caption=false,font=footnotesize]{subfig}
\usepackage{graphbox} 

\usepackage{algpseudocode, algorithm, algorithmicx}
\newcommand*\Let[2]{\State #1 $\gets$ #2}
\algrenewcommand\algorithmicrequire{\textbf{Given:}}
\algrenewcommand\algorithmicensure{\textbf{Return:}}

\usepackage[english]{babel}
\usepackage{bm}
\usepackage{hyperref}

\usepackage{fancyhdr}
\fancypagestyle{titlepage}{
\fancyhead{\begin{minipage}{\textwidth}
\footnotesize\copyright 2021 IEEE. Personal use of this material is permitted. Permission from IEEE must be obtained for all other uses, in any current or future media, including reprinting/republishing this material for advertising or promotional purposes, creating new collective works, for resale or redistribution to servers or lists, or reuse of any copyrighted component of this work in other works.\newline
This article is published in IEEE Xplore (\href{https://ieeexplore.ieee.org/document/9661310}{DOI: 10.1109/TSP.2021.3137599})\vspace*{-10pt}
\end{minipage}
}
}

\newtheorem{theorem}{Theorem}
\newtheorem{corollary}[theorem]{Corollary}
\newtheorem{lemma}[theorem]{Lemma}
\newtheorem{proposition}[theorem]{Proposition}
\newtheorem{definition}[theorem]{Definition}
\newtheorem{remark}[theorem]{Remark}

\renewcommand{\a}{\mathbf{a}}
\newcommand{\A}{\mathbf{A}}

\newcommand{\B}{\mathbf{B}}
\newcommand{\C}{\mathbf{C}}

\newcommand{\D}{\mathbf{D}}

\newcommand{\Pb}{\mathbf{P}}

\newcommand{\Rc}{\mathcal{R}}
\renewcommand{\S}{\mathbf{S}}

\newcommand{\x}{\mathbf{x}}
\newcommand{\X}{\mathbf{X}}

\newcommand{\Xzo}{\mathbf{X}_{01}}
\newcommand{\Xzov}{\vv{\mathbf{X}}_{01}}
\newcommand{\y}{\mathbf{y}}
\newcommand{\Y}{\mathbf{Y}}
\newcommand{\z}{\mathbf{z}}
\newcommand{\Z}{\mathbf{Z}}

\newcommand{\Zzov}{\vv{\mathbf{Z}}_{01}}

\newcommand{\eps}{\varepsilon}

\newcommand{\LAmbda}{{\boldsymbol{\lambda}}}
\newcommand{\THeta}{\bm{\vartheta}}

\newcommand{\xxi}{\bm{\xi}}
\newcommand{\0}{\bm{0}}

\newcommand{\id}{\mathbf{Id}}

\renewcommand{\vec}{\mathrm{vec}}

\newcommand{\sign}{\mathrm{sign}}
\newcommand{\supp}{\mathrm{supp}}
\newcommand{\corr}{\mathrm{Corr}}

\newcommand{\argmin}{\operatornamewithlimits{argmin}}
\newcommand{\R}{\mathbb{R}}

\newcommand{\curly}[1]{\left\{ #1 \right\}}
\newcommand{\abs}[1]{\left| #1 \right|}
\newcommand{\inner}[1]{\left\langle #1 \right\rangle}
\newcommand{\norm}[1]{\left\| #1 \right\|}

\begin{document}
	
\title{Structural Sparsity in Multiple Measurements}
\author{F. Bo\ss{}mann, S. Krause-Solberg, J. Maly, and N. Sissouno
\thanks{FB acknowledges funding by the National Science Foundation of China (NSFC), project number 42004109 (Seismic inpainting using wave models on scattered data).
SKS acknowledges support by the Helmholtz Imaging Platform (HIP), a platform of the Helmholtz Incubator on Information and Data Science.
JM  acknowledges funding  by the Deutsche  Forschungsgemeinschaft  (DFG,  German  Research  Foundation) through the  project CoCoMIMO funded  within  the  priority  program SPP 1798 \emph{Compressed  Sensing in  Information Processing} (COSIP).
The figures in Section \ref{sec:meteor} were made using \textit{M\_Map}, a mapping package for MATLAB by R. Pawlowicz (online available at \url{https://www.eoas.ubc.ca/~rich/map.html}).}%

\thanks{F. Bo\ss{}mann, Harbin Institute of Technology, Department of Mathematics, f.bossmann@hit.edu.cn}
\thanks{S. Krause-Solberg, Deutsches Elektronen-Synchrotron DESY, Notkestra\ss{}e 85, D-22607 Hamburg, sara.krause-solberg@desy.de}%
\thanks{J. Maly, KU Eichst\"att, Ostenstra\ss{}e 26, D-85072 Eichst\"att, \mbox{johannes.maly@ku.de}}%
\thanks{N. Sissouno, Technical University of Munich, Faculty of Mathematics, Boltzmannstra\ss{}e 3,
D-85748 Garching, sissouno@ma.tum.de}%
}

\maketitle
\thispagestyle{titlepage}
\begin{abstract}
	We propose a novel sparsity model for distributed compressed sensing in the multiple measurement vectors (MMV) setting. Our model extends the concept of row-sparsity to allow more general types of structured sparsity arising in a variety of applications like, e.g., seismic exploration and non-destructive testing. To reconstruct structured data from observed measurements, we derive a non-convex but well-conditioned LASSO-type functional. By exploiting the convex-concave geometry of the functional, we design a projected gradient descent algorithm and show its effectiveness in extensive numerical simulations, both on toy and real data.
\end{abstract}

\IEEEpeerreviewmaketitle


\section{Introduction}
Starting with the seminal works \cite{candes2006robust,candes2006near,donoho2006compressed} a rich theory on signal reconstruction from seemingly incomplete information has evolved under the name \emph{compressed sensing} in the past two decades, cf.\ \cite{foucart2013mathematical} and references therein. 

The core idea is to use the intrinsic structure of a high-dimensional signal $\x \in \R^N$ to allow reconstruction from $m \ll N$ linear measurements
\begin{align} \label{eq:CS}
    \y = \A \x,
\end{align}
where $\A \in \R^{m\times N}$ models the measurement process and $\y \in \R^m$ is called the \emph{measurement vector}. In what follows we will call the columns $\a_j$ of $\A$ \emph{atoms}. One particular instance of intrinsic structure is \emph{sparsity}: the signal $\x$ is called $s$-sparse if $\abs{\supp(\x)} \le s$, where $\supp(\x) = \curly{i \colon x_i \neq 0}$ denotes the \emph{support} of $\x$. We will call an atom $\a_j$ \emph{activated} if $j \in \supp(\x)$. If $\A$ is well-designed, $m \approx s \log(\tfrac{N}{s})$ measurements suffice to guarantee stable and robust recovery of all $s$-sparse $\x$ from $\y$ by polynomial time algorithms \cite{foucart2013mathematical}. This suggests that the number of necessary measurements mainly depends on the intrinsic information encoded in $\x$.\\
Despite its simplicity the model in \eqref{eq:CS} encompasses many real-world measurement set-ups. For instance, in seismic exploration, a key challenge is to reconstruct earth layers from few linear measurements \cite{Bossmann15,Bossmann16}. In this case, the vector $\x$ is a discretized vertical slice through the ground where each entry represents the seismic reflectivity. To reconstruct the earth layers, a synthetic seismic impulse is produced and its reflections are measured at different positions on the surface. This measurement process can be modeled by a convolution of $\x$ with the seismic impulse such that $\A$ is the corresponding convolution matrix. Since the reflectivity is low whenever the material is mostly homogeneous and high at material boundaries, the vector $\x$ can be assumed to be sparse and its non-zero entries indicate the earth layer boundaries. A similar model applies to ultrasonic non-destructive testing where an ultrasonic impulse is sent into an object and defects inside the material are reconstructed from the reflections of this impulse \cite{Bossmann12}. Other possible applications are face and speech recognition \cite{Katz87,Wright09}, magnetic resonance imaging \cite{Lustig07}, or computer tomography \cite{Sidky08}. For an overview also see \cite{Starck10,Rao98} and the references therein.\\
In all applications mentioned above, we can assume structure not only in one direction of space but in multiple dimensions, meaning that measurements at different locations/times $t_1 < ... < t_L$ correspond to different ground-truth signals $\x_1,...,\x_L \in \R^N$ whose support structure is related. When thinking of waves travelling through the ground, the positions of non-zero entries of consecutive $\x_l$ can only differ up to a certain number determined by properties of the surrounding material and fineness of the discretization.\\
If several signals $\x_1,...,\x_L \in \R^N$ are measured by the same process $\A$, the model in \eqref{eq:CS} becomes
\begin{align} \label{eq:MultipleMeas} 
	\Y = \A\X,
\end{align} 
for $\X=(\x_1,\ldots,\x_L) \in \R^{N\times L}$ and $\Y=(\y_1,\ldots,\y_L) \in \R^{m \times L}$. In this setting, also known as \emph{multiple measurement vectors (MMV)}, the necessary number of measurements may be reduced by exploiting joint structure in $\X$, for instance, row sparsity (all $\x_l$ share a common support). This has already been done in applications like MRI \cite{wu:2014} and MIMO communications \cite{rao:2014}. In the case of row-sparsity, reconstruction is usually performed via
\begin{align}\label{eq:row-0}
    \min_{\X \in \R^{N\times L}} \| \A\X-\Y \|_F^2 + \lambda\|\X\|_{\text{row}-0}\ ,
\end{align}
where $\| \cdot \|_F$ denotes the Frobenius-norm, $\|\cdot\|_{\text{row-0}}$ denotes by abuse of notation the number of non-zero rows, and $\lambda$ is a tunable parameter. Although \eqref{eq:row-0} is NP-hard in general \cite{foucart2013mathematical}, solutions can be well approximated via greedy algorithms or convex relaxation.

\begin{figure}
    \centering
    \includegraphics[align=c,width=2.8cm]{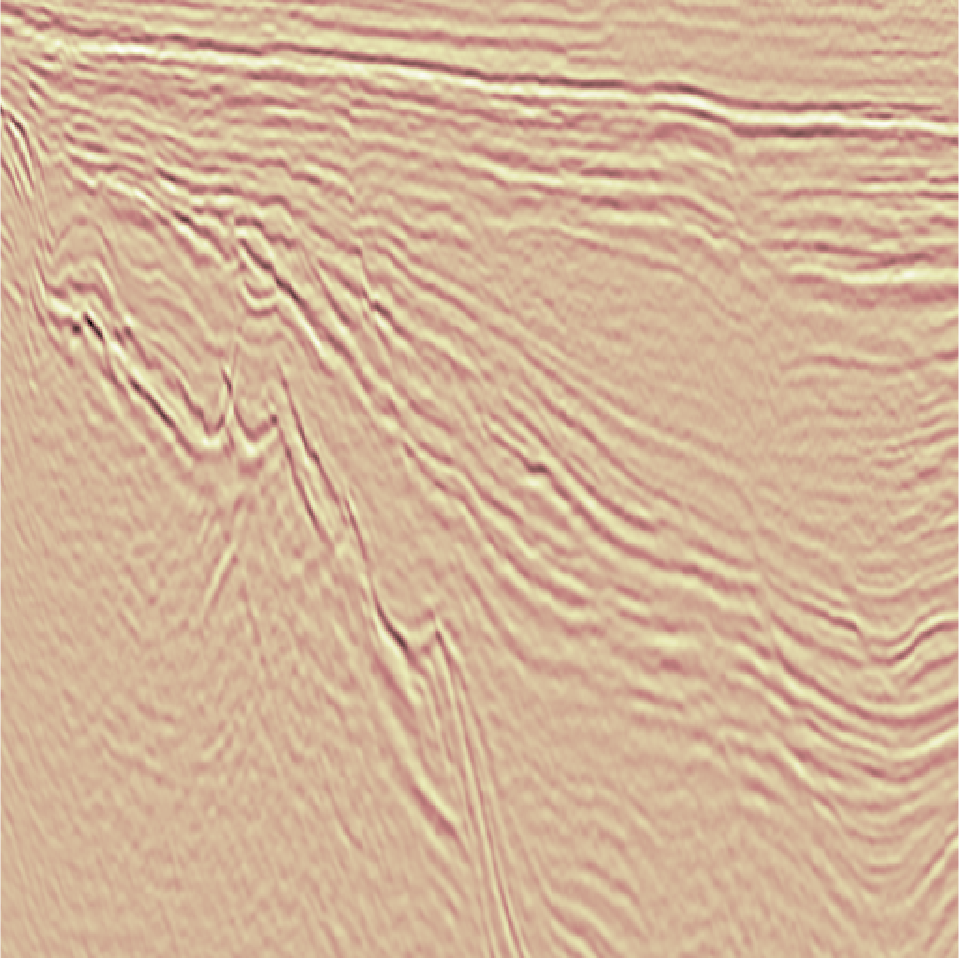}
    \includegraphics[align=c,width=2.8cm]{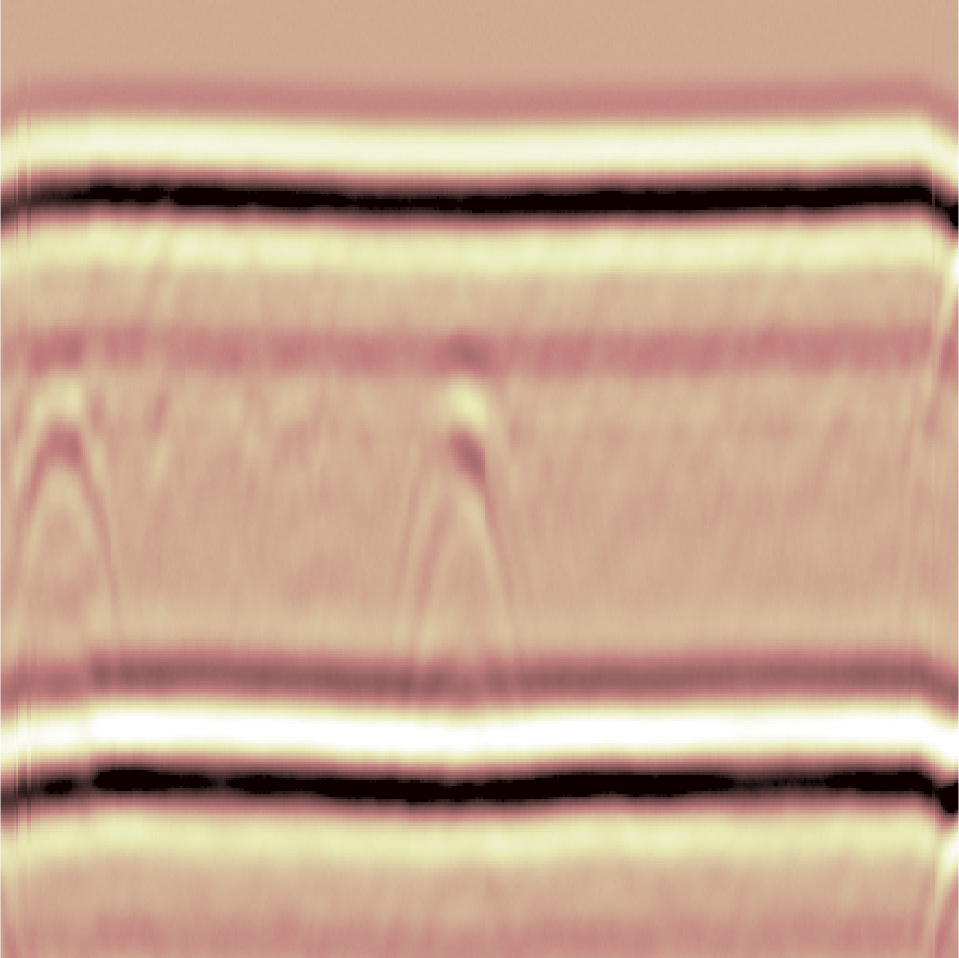}
    \includegraphics[align=c]{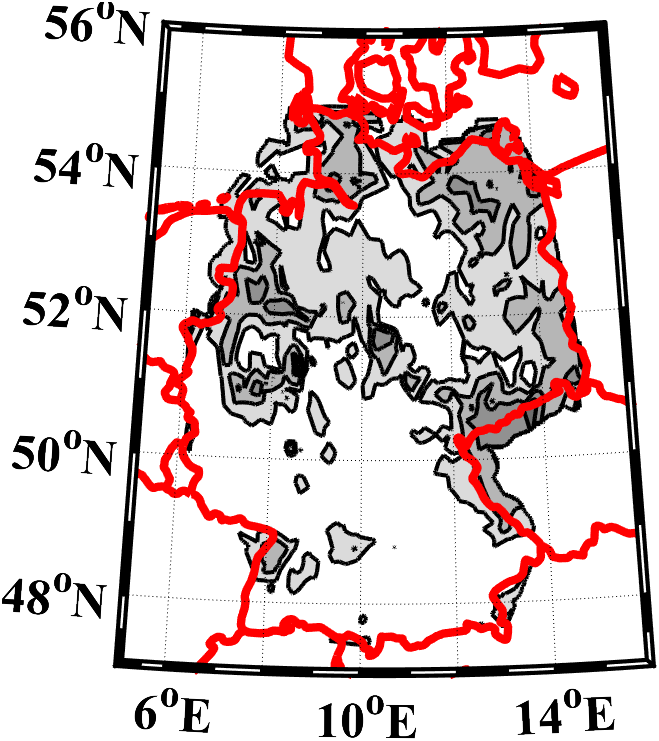}
    \caption{Data from different applications exhibiting clear structure but being neither row- nor block-sparse: Seismic exploration (left), ultrasonic non-destructive testing (middle), and meteorology (right).}
    \label{fig:seismicExample}
\end{figure}

However, row-sparsity and other established structural models (column sparsity, block sparsity) make restrictive assumptions on the concrete structure of $\X$. As a simple thought experiment, consider $\X$ to be the identity matrix. Then $\X$ is neither row- nor block-sparse but clearly exhibits a simple structure, namely a diagonal line. In fact, the established models are too restrictive for many applications. For example, the support and thus the structure may change over time as it is the case in real-time dynamic MRI \cite{Vaswani2010}, dynamic PET \cite{Heins2015}, and wireless communication \cite{Lian2019}. In machine learning as well as in quantile and logistic regression modeling more sophisticated structure models may be required \cite{Lee2014,Huang2011}. In this work, we consider the three applications seismic exploration, ultrasonic non-destructive testing, and meteorology. Here, the data is gathered at different locations where the support might change over the spatial dimension. For instance, the above mentioned earth layers in seismic exploration do not follow straight horizontal lines and thus are not even close to row sparse. Material defects that have to be reconstructed in non-destructive testing can have many forms, most of which are neither row nor block sparse. In Fig. \ref{fig:seismicExample} we show some exemplary measurements that do not fit into established structural models.

 Moreover, in some applications there does not exist a left-to-right order of the measurement vectors $\y_l$  (for instance, when the measurements are taken from scattered locations) and also the order of atoms $\a_j$ may not be clear. In this case, we want to reconstruct the solution independent of column permutations in $\A$ and $\Y$. If $\X$ is the structural sparse solution of $\A\X=\Y$, then $\tilde{\X}=\Pb_\Y\X\Pb_\A$ should be the structural sparse solution of $(\Pb_\A \A)\tilde{\X}=\Pb_\Y \Y$, where $\Pb_\A$ and $\Pb_\Y$ are permutation matrices. This excludes all order-dependent approaches to define structural sparsity since such definitions are not invariant under permutations.

\subsection{Contribution} 

In this work, we introduce a sparsity model that can capture a wide range of practically relevant structures, comes with efficient optimization, and allows to learn the structures in an intuitive way from only the measurements and additional knowledge of the concrete application. To be more precise, our model is a special case of group sparsity for multiple measurement recovery problems and encompasses established concepts like row- and block-sparsity. The novel ingredient is that the structural support constraints are encoded in a matrix $\C$ which allows efficient processing. We introduce a non-convex regularizer enforcing the structures encoded in $\C$ and discuss possible relaxations of the regularizer. Based on theoretical insights into the optimization landscape of our regularizer, we suggest a projected gradient descent to minimize the related LASSO-type formulation. Finally, we provide a simple heuristic to determine $\C$ for concrete  applications under sole knowledge of the measurement process and measurements. We validate efficacy of both the parameter heuristic and the regularizer in extensive numerical simulations on real data.

\subsection{Related Work}

There exist several approaches to solve \eqref{eq:MultipleMeas} by assuming different sparsity models for $\X$. Most of them adapt methods that were originally designed to solve \eqref{eq:CS}. In \cite{Cotter05} an extension of the algorithms Matching Pursuit and the FOCal Underdetermined System Solver (FOCUSS) are presented. Bayesian methods are considered in \cite{Zhang10,Ziniel13}. In \cite{Tropp06a,Tropp06b} the authors introduce greedy pursuits and convex relaxations for the MMV problem. Theoretical results have been shown, e.g., in \cite{Chen06}. All these methods enforce row-sparsity in the reconstructed solution. In \cite{Baron06} two joint sparsity models (JSM) for compressed sensing are introduced. JSM-1 considers solutions where all columns can be written as the sum of a common sparse component that is equal for each column and another unique sparse vector. This is equivalent to assuming $\X=\X_1+\X_2$ where $\X_1$ is row-sparse and $\X_2$ is sparse (without any correlation between different columns). JSM-2 is a slight relaxation of row-sparsity and allows small support changes over the columns. Yet another approach is presented in \cite{Lu12,Zheng11}: correlated measurements are assumed to have sparse approximations that are close in the Euclidean distance. This idea is related to dynamic compressed sensing \cite{Angelosante09,Ziniel13a} where neighboring columns are assumed to have similar support. In both cases, the support is allowed to change slowly over different data vectors. Nevertheless, the above methods share quite restrictive support assumptions based on geometrical features and hence cannot reconstruct simple features in the solution that do not match those strict assumptions.\\
A more general approach is the so-called group sparsity model \cite{baraniuk2010model,huang2010benefit} in which a set of groups $\mathcal{G}$ is defined whose elements $G \in \mathcal{G}$ encode support sets. The matrix $\X$ is called $s$-group-sparse if $\supp(\X)$ is a subset of a union of at most $s$ groups in $\mathcal{G}$. Whereas this model is able to encode all possible structural constraints on $\X$, its generality comes with a price. First, straight-forward adaption of compressed sensing algorithms is only possible under knowledge of the set $\mathcal{G}$, cf.\ \cite{baraniuk2010model,huang2010benefit} and subsequent literature. Second, even if $\mathcal{G}$ is known, its cardinality might grow exponentially in the ambient dimension which considerably increases the computational complexity of established procedures. Third, if $\mathcal{G}$ is unknown, learning algorithms have to either rely on clustering of entries \cite{yu2012bayesian,zhang2013extension} leading to block-sparse-like models, or on available training data $(\X_i,\Y_i)$, for $i \in [n]$, from which $\mathcal{G}$ can be learned \cite{peleg2012exploiting}. The latter work suggests an alternating approach to learn both signal and structure without initial data; nevertheless, it does not provide as simple means to incorporate expert knowledge on concrete applications as our heuristic for determining the structure encoding matrix $\C$.
Let us finally mention that, for general $\mathcal{G}$, convex regularizers can only be computed theoretically, cf.\ the concept of atomic norms in \cite{chandrasekaran2010convex}. \\

\begin{remark}
   In this work we concentrate on understanding and representing the intrinsic sparsity structure of $\X$ for a fixed measurement process/dictionary $\A$. Other lines of work discuss how to learn a proper dictionary for $\X$ (such that $\X$ becomes sparse in a classical sense) or how to adapt an existing one by slight perturbation to improve reconstruction performance \cite{aharon2006k,zhu2011sparsity}. They, however, do not allow to represent general structural dependencies between active entries of $\X$. Combining those approaches with our generalized sparsity model is an interesting topic for future work.
\end{remark}

\subsection{Notation and Outline} 

We denote matrices by bold capital letters, vectors by bold lowercase letters, and scalars by regular letters.
The only exceptions are vectors that we get by the vectorization $\vv{\Z} := \vec(\Z)$ of matrices $\Z$. The inversion (reshape) of the vectorization is denoted by $\vec^{-1}$.
As already mentioned in the introduction, the columns of a matrix $\Z\in\R^{m\times n}$ are denoted by $\z_l$ for
$l\in [n]$, where $[n]:=\{1,\dots,n\}$ is used to abbreviate index sets. The identity matrix and the matrix of ones are written as $\id$ and $\mathbb{1}$, respectively.
For the set of non-negative real numbers we use the notation $\R_+$.\\
We denote the support matrix of $\Z$ by $\Z_{01}$, i.e., $\Z_{01}\in\{0,1\}^{m\times n}$ is the matrix with entries $|\sign(Z_{j,l})|$ for
$j\in [m]$ and $l\in [n]$. If applied to matrices or vectors, the $\sign$-function as well as the absolute value $|\cdot|$ act component-wise.

Besides the standard matrix multiplication we use the Kronecker product $\otimes$ defined by $\A\otimes \B:=(A_{j,l}\B)_{j,l}\in\R^{mp\times nq}$, for matrices $\A\in\R^{m\times n}$ and $\B\in\R^{p\times q}$. The Frobenius norm of a matrix $\Z$ is $\|\Z\|_F$, while $\|\z\|_2$ denotes the Euclidian norm of a vector $\z$.

The outline of the paper is as follows. In Section~\ref{sec:StructuralSparsity}, we introduce a general model for structural sparsity and discuss possible reconstruction approaches of such structured signals. In particular, we derive a relaxed, still non-convex functional, whose minimizers provide good approximations to $\X$, and explore the specific geometry of the functional. Building upon these insights, we describe in Section~\ref{sec:Optimization} a projected gradient descent procedure to efficiently solve the program. Finally, in Section~\ref{sec:Numerics}, we empirically validate our model on toy scenarios and real (seismic/ultrasonic/meteorological) data.


\section{Structural Sparsity} \label{sec:StructuralSparsity}
We begin by deriving a general model for structural sparsity. After discussing its relation to established structural models like row- or column-sparsity, we provide a corresponding NP-hard optimization problem to reconstruct structured signals from compressive measurements. To solve the in general intractable problem, we suggest a non-convex relaxation of particular convex-concave shape.

\subsection{The Basic Model} \label{sec:BasicModel}

In order to develop a notion of structural sparsity that is capable of describing signals like the ones in Fig. \ref{fig:seismicExample}, we first have to understand the underlying abstract idea of row-sparsity and related concepts. A matrix $\X$ is $s$-row-sparse if it has at most $s$ non-zero rows, i.e., if there are up to $s$ matrices $\X_k$ with exactly one non-zero row such that $\X = \X_1 + \cdots + \X_s$. We could say that the matrices $\X_k$ describe the \emph{elementary structures} of row-sparsity. By changing the elementary structures, one obviously recovers various established concepts like sparsity ($\X_k$ are matrices with exactly one non-zero entry) and block-sparsity ($\X_k$ are matrices with exactly one non-zero block). This elementary idea is the corner stone of group sparsity \cite{baraniuk2010model,huang2010benefit}.

Building upon the same intuition, we wish to describe elementary structures in a practical way that lends itself to efficient computation. To this end, given two non-zero entries $X_{j,l},X_{j',l'}$ of a structured signal matrix $\X \in \R^{N \times L}$, let $C_{(j,l),(j',l')}\in \curly{0,1}$ indicate whether the two entries can belong to a single elementary structure of $\X$. To be more precise, $C_{(j,l),(j',l')} = 0$ if they can belong to the same structure, i.e., there exists (at least) one elementary structure $\X_k$ whose entries $(j,l)$ and $(j',l')$ are non-zero. If such a structure does not exist, we set $C_{(j,l),(j',l')} = 1$. Then,
\begin{align}\label{eq:pen1}
    \sum\limits_{\substack{j,j',l,l' \colon \\ X_{j,l},X_{j',l'}\neq0}} C_{(j,l),(j',l')}
    = 0
\end{align}
whenever $\X$ itself is an elementary structure. Moreover, the value of (\ref{eq:pen1}) increases the more $\X$ differs from an elementary structure. Let us clarify this by a simple example: the choice $C_{(j,l),(j',l')}=0$ for $j=j'$ and $1$ otherwise would characterize the basic units of row-sparsity as \eqref{eq:pen1} is $0$ if and only if $\X$ has at most one non-zero row. 
Using the vectorization of the support matrix $\Xzo \in \{0,1\}^{N\times L}$ of $\X\in\R^{N\times L}$ we can rewrite (\ref{eq:pen1}) as
\begin{align} \label{eq:pen2}
    \Xzov^T\C \Xzov = 0,
\end{align}
where $\C \in \R^{NL \times NL}$ has entries $C_{(j,l),(j',l')}$, for $j,j' \in [N]$ and $l,l' \in [L]$. We may now define the set of $\C$-structured $s$-sparse signals
\begin{align} \label{eq:Scs}
	\mathbb{S}_\C^s &= \left\{\Z \in \R^{N\times L}\colon
	\begin{array}{cc}
	\Z = \sum_{k = 1}^{s} \Z_k,\\
	(\vv{\Z}_k)_{01}^T \C (\vv{\Z}_k)_{01} = 0, \; \forall k \in [s]
	\end{array}
	\right\}.
\end{align}
\begin{remark} \label{rem:Models}
	The model defined in \eqref{eq:Scs} is quite general and covers several well-known special cases (in addition to row-sparsity mentioned above). Choosing $\C = \mathbb{1} - \id$, the set $\mathcal{S}_\C^s$ describes the set of $s$-sparse vectors. Choosing $\C$ such that
	\begin{align*}
		C_{(j,l)(j',l')} = \begin{cases}
		0 & |j - j'| \le a \text{ and } |l - l'| \le b, \\
		1 & \text{else,}
		\end{cases}
	\end{align*}
	we recover the set of block-sparse matrices with blocks of size $a\times b$ (cf.\ \cite{baraniuk2010model}).\\
	Note that the diagonal entries of $\C$ are $0$ independent of the concrete model as $\C$ shall only characterize relations between different entries of $\X$.
\end{remark}
Building upon \eqref{eq:Scs} we can define the $\C$-structured $\ell_0$-norm
\begin{align} \label{eq:lC0}
	\| \Z \|_{\C,0} = \min \{s \ge 0 \colon \Z \in \mathbb{S}^s_\C \},
\end{align}
which is actually not a norm but abuse of notation. Minimizing the $\ell_{\C,0}$-norm constrained to correct measurements, i.e.,
\begin{align} \label{eq:lC0min}
	\min_{\Z \in \R^{N\times L}} \| \Z \|_{\C,0}, \quad \text{subject to } \A\Z = \Y,
\end{align}
then extends $\ell_0$-minimization to the $\C$-structured $s$-sparse case. The program in \eqref{eq:lC0min} inherits NP-hardness from classical sparse recovery such that it is undesirable to solve \eqref{eq:lC0min} directly. In fact, even computing \eqref{eq:lC0} is NP-hard in general. Note, however, that \eqref{eq:pen2} itself might suffice as regularizer since its magnitude increases if the number of elementary structures in $\X$ increases. This handwavy argument is substantiated by the observation that, under mild assumptions on the elementary structures encoded in $\C$, there is an equivalence relation between \eqref{eq:pen2} and \eqref{eq:lC0} as stated in the following proposition.

\begin{proposition} \label{thm:Equivalence}
	Assume that, for $s_\text{col} \in [N]$, the sparsity model characterized by $\C$ satisfies
	\begin{align} \label{eq:Assumption}
	    \| \X \|_{\C,0} = 1 \implies \| \x_l \|_0 \le s_\text{col}, \; \forall l \in [L].
	\end{align}
	Then, for $\X \neq \0$, we have
	\begin{align*}
		\frac{1}{2} \| \X \|_{\C,0}^2 \le \Xzov^T \C \Xzov + 1 \le s_\text{col}^2 L^2 \| \X \|_{\C,0}^2.
	\end{align*}
\end{proposition} 
\begin{IEEEproof}
    For $\| \X \|_{\C,0} = 1$, one has that $\Xzov^T \C \Xzov = 0$ by \eqref{eq:Scs}. Let us now assume $\| \X \|_{\C,0} = s \ge 2$. Then we can write $\Xzo = \sum_{k=1}^{s} \X_k$ where each $\X_k$ contains only ones and zeros. We hence get
    \begin{align*} 
    	\Xzov^T \C \Xzov = \left( \sum_{k=1}^{s} \vv{\X}_k \right) \C \left( \sum_{k=1}^{s} \vv{\X}_k \right) = \sum_{j \neq k} \vv{\X_k}^T \C \vv{\X_j},
    \end{align*} 
    since $\vv{\X_k}^T \C \vv{\X_k} = 0$, for all $k \in [s]$. By assumption, the matrices $\X_k$ have at most $s_\text{col} L$ non-zero entries such that
    \begin{align} \label{eq:proofInequality}
    	1 \le \vv{\X_k}^T \C \vv{\X_j} \le s_\text{col}^2 L^2,
    \end{align} 
    where the lower bound is a consequence of the minimal decomposition required in \eqref{eq:lC0}. We conclude that
    \begin{align*}
    	s(s-1) \le \sum_{j \neq k} \vv{\X_k}^T \C \vv{\X_j} \le s_\text{col}^2 L^2 s(s-1),
    \end{align*}
    which yields the claim.
\end{IEEEproof}
\begin{remark} \label{rem:Assumption}
    Assumption \eqref{eq:Assumption} requires the elementary structures described by $\C$ to have $s_\text{col}$-sparse columns which holds for $s_\text{col} = 1$ when working with generalizations of row-sparsity. This can be easily seen in Fig. \ref{fig:seismicExample}. Let us emphasize that the construction of $\C$ presented in Appendix \ref{sec:C} satisfies \eqref{eq:Assumption} if the appearing parameters $\alpha$ and $\beta$ are suitably chosen.
\end{remark}

If we have information on the structure of $\X$ in form of $\C$ and the measurements $\A$ are injective when restricted to $\mathbb{S}_\C^s$, Proposition \ref{thm:Equivalence} suggests to approximate $\X$ from $\Y$ by solving
\begin{align} \label{eq:BP}
	\min_{\Z \in \R^{N \times L}} \Zzov^T \C \Zzov, \quad \text{subject to } \A\Z = \Y.
\end{align}
This raises two questions: how can one obtain a suitable structure matrix $\C$ in general and how can the NP-hard optimization in \eqref{eq:BP} be solved? We refer the interested reader to Appendix \ref{sec:C} for a simple heuristic to construct $\C$ from $\A$ and $\Y$, and now address the second question.

\subsection{Ways of Relaxation}
\label{sec:WaysOfRelaxation}

The program in \eqref{eq:BP} poses two difficulties. First, the matrix $\C$ is not necessarily positive semi-definite and, second, the support vector $\Zzov$ depends in a non-continuous way on $\Z$. 
To circumvent the latter, we define the regularizer 
\begin{align*}
    \Rc_{\C} \colon \R^{N \times L} \rightarrow \R, \quad \Rc_{\C} (\Z) = \vv{\Z}^T \C \vv{\Z},
\end{align*}
and rewrite \eqref{eq:BP} as the binary program
\begin{align}\label{BinaryProgram}
    \min_{\substack{\tilde{\Z} \in \{ 0,1 \}^{N \times L}, \\ \sign(|\Z|) = \tilde{\Z}, \\ \A\Z = \Y}} \Rc_{\C} (\tilde{\Z}),
\end{align}
where $\sign$ as well as $|\cdot|$ act component-wise on matrices. As \eqref{BinaryProgram} illustrates, the non-continuous/non-convex dependence of $\Zzov$ on $\Z$ lies in the relation between $\Z$ and the auxiliary variable $\tilde \Z$. Replacing the $\sign$-function by identity (interpreting identity as a convex relaxation of $\sign$) leads to
\begin{align} \label{eq:BPrelaxed}
	\min_{\Z \in \R^{N \times L}} \Rc_{\C} (|\Z|), \quad \text{subject to } \A\Z = \Y.
\end{align}
Since \eqref{eq:BPrelaxed} does not incorporate noise on the measurements, we replace it with the more robust formulation
\begin{align} \label{eq:LASSO}
	\min_{\Z \in \R^{N \times L}} \| \A\Z - \Y \|_F^2 + \lambda \Rc_{\C} (|\Z|),
\end{align}
where $\lambda > 0$ is a regularization parameter. It is well known from related programs that solutions of \eqref{eq:LASSO} solve a robust version of \eqref{eq:BPrelaxed} while the magnitude of $\lambda$ balances robustness and accuracy of the reconstruction. We detail this in the following lemma. The proof is similar to \cite[Proposition 3.2]{foucart2013mathematical} and thus omitted.
\begin{lemma} \label{lem:Lasso}
	If $\X_\lambda$ minimizes \eqref{eq:LASSO} with parameter $\lambda > 0$, then $\X_\lambda$ minimizes
	\begin{align*}
		\min_{\Z \in \R^{N \times L}} \Rc_{\C} (|\Z|), \quad \text{subject to } \| \A\Z - \Y \|_F \le \eta_\lambda,
	\end{align*}
	where $\eta_\lambda = \| \A\X_\lambda - \Y \|_F$.
\end{lemma}

Although \eqref{eq:LASSO} is a non-convex problem, the regularizer $\Rc_\C$ exhibits some beneficial geometrical properties.

\begin{proposition} \label{thm:Geometry}
    For any $\X,\D \in \R^{N\times L}$ the following holds. Let $\S = \sign(\X) \in \R^{N\times L}$ and define $\R_\S^{N\times L}$ as the orthant of $\R^{N\times L}$ in which $\X$ lies. Then,
    \begin{align} \label{eq:hFunction}
    \begin{split}
        f \colon \curly{t \in \R \colon \X + t\D\in\R_\S^{N\times L} } \rightarrow \R,\\
        f(t) = \Rc_\C (\abs{\X + t\D})
    \end{split}
    \end{align}
    is a convex or concave function. In particular, if $\D \in \R_\S^{N\times L}$ or $\D \in -\R_\S^{N\times L}$, then the function in \eqref{eq:hFunction} is convex.
\end{proposition}
\begin{IEEEproof}
    Recall that all entries of $\C$ are non-negative and that $\C$ is symmetric by definition. Obviously,
    \begin{align*}
        g(t) = \Rc_\C(\X + t\D) = t^2 \vv{\D}^T \C \vv{\D} + 2t \vv{\X}^T \C \vv{\D} + \vv{\X}^T \C \vv{\X}
    \end{align*}{}
    is a quadratic functional and thus either convex or concave. If $\D \in \pm \R_+^{N\times L}$, the functional $g$ is convex as $\vv{\D}^T \C \vv{\D} \ge 0$. By restricting $f$ such that $\X+t\vv{\D}\in\R_\S^{N\times L}$, we get that
    \begin{align*}
        f(t) &= \Rc_\C(\S \odot (\X + t\D)) \\&= t^2 \vv{\D}^T \tilde{\S}^T \C \tilde{\S} \vv{\D} + 2t \vv{\X}^T \tilde{\S}^T \C \tilde{\S} \vv{\D} + \vv{\X}^T \tilde{\S}^T \C \tilde{\S} \vv{\X},
    \end{align*}{}
    where $\odot$ denotes the Hadamard product and $\tilde{\S} \in \R^{NL\times NL}$ the diagonal matrix with $\vv{\S}$ on its diagonal. Obviously, the restriction of $f$ equals $g$ where $\C$ is replaced by $\tilde{\S} \C \tilde{\S}$ which gives the first claim. The second claim follows since $\tilde{\S} \vv{\D} \in \pm \R_+^{NL}$, for $\D \in \pm \R_\S^{N\times L}$.
\end{IEEEproof}
Proposition \ref{thm:Geometry} states that if restricted to single orthants, $\Rc_\C(\abs{\cdot})$ behaves well along rays. In particular, the function $\Rc_\C(\abs{\cdot})$ is convex along all rays passing through the origin.

\begin{corollary} \label{cor:Geometry}
    For any $\D \in \R^{N\times L}$, the function
    \begin{align*}
        t \mapsto \Rc_\C (\abs{t\D})
    \end{align*}{}
    is convex.
\end{corollary}{}

An important consequence of Proposition \ref{thm:Geometry} is that if one had oracle knowledge on the orthant of $\X$, the program in \eqref{eq:LASSO} could be restricted accordingly and would become better conditioned.\\
The naive approach thus would be to pick any initialization sharing the same sign with $\X$ and then restricting \eqref{eq:LASSO} to the corresponding orthant. Unfortunately, the possibly most common initialization, the back-projection of $\Y$ by the pseudo-inverse $\A^\dagger$ of $\A$
\begin{align} \label{eq:PseudoInverse}
	\X_0 = \A^\dagger \Y = \argmin_{\Z \in \R^{N \times L}} \| \Z \|_F, \quad \text{subject to } \A\Z = \Y,
\end{align} 
in general does not have this property as the following theorem shows\footnote{Although the initialization in \eqref{eq:PseudoInverse} does not always provide the necessary orthant information, we mention that it often succeeds in numerical experiments.}.

\begin{theorem} \label{thm:Counterexample}
	For any (at least $2$-sparse) vector $\x_0 \in \R^N$ (here $L = 1$), there exists $\A \in \R^{N\times N}$ and $\y \in \R^N$ such that $\x_0$ fulfills \eqref{eq:PseudoInverse} and there exists a $2$-sparse vector $\x^{\text{sp}}$ which solves the linear system, i.e., $\A\x^{\text{sp}} = \y$, but lies in another orthant.
\end{theorem}
\begin{IEEEproof}
    Without loss of generality let $\x_0 \geq 0$ (entry-wise) and assume that its two first entries are non-zero. Define  $\x^{\text{sp}} = (-a,b,0,0,\dots)$ with $a,b>0$ and $\THeta = \x_0 - \x^{\text{sp}}$. Choose $a,b > 0$ such that
    \begin{align*}
       0 = \langle \THeta,\x_0 \rangle = \| \x_0 \|_2^2 + a (x_0)_1 - b (x_0)_2,
    \end{align*}
    which is equivalent to
    \begin{align} \label{eq:b}
       b = \frac{\| \x_0 \|_2^2+ a (x_0)_1}{(x_0)_2}.
    \end{align}
    Now, let $\A \in \R^{N\times N}$ be any matrix with $\mathrm{span}( \THeta ) = \ker(\A)$ and define $\y = \A \x_0$. Then $\x_0$ is perpendicular to the kernel and thus the minimum norm solution of \eqref{eq:PseudoInverse}. Furthermore, $\A\x^{\text{sp}} = \y$ but $\x^{\text{sp}}$ lies in a different orthant than $\x_0$.
\end{IEEEproof}

Instead of searching for suitable alternative initialization procedures, we propose to modify the optimization problem. Indeed, we can embed \eqref{eq:MultipleMeas} into higher dimensional spaces and work with the augmented linear system
\begin{align} \label{eq:Extended}
	\A^\pm \Z^\pm = \Y
\end{align}
where $\A^\pm = (\A, -\A) \in \R^{m \times 2N}$ and $\Z^\pm \in \R^{2N \times L}$. Obviously, the original solution $\X$ solves \eqref{eq:Extended} by defining
\begin{align}
	\X^\pm = \begin{pmatrix}
	\X \\
	\0
	\end{pmatrix}.
\end{align}
More importantly, if we define $\X_+ \in \R_+^{N \times L}$ and $\X_- \in \R_+^{N \times L}$ as positive and negative part of $\X$ (containing only the positive/negative entries of $\X$ in absolute value and setting the rest to zero such that $\X = \X_+ - \X_-$), we have that the matrix
\begin{align*}
	\X^\pm_{\text{pos}} = \begin{pmatrix}
	\X_+ \\
	\X_-
	\end{pmatrix}
\end{align*}
solves \eqref{eq:Extended}, shares the structural complexity of $\X$, and lies within the positive orthant. To summarize the last lines, increasing the dimension of the linear system only by a factor of two, we can guarantee that a structured solution (from which the original $\X$ is directly recovered) may be found by applying an appropriate solver to \eqref{eq:LASSO} restricted to the positive orthant. From now on, we hence assume without loss of generality that $\X$ itself lies within the positive orthant and consider the restricted program
\begin{align} \label{eq:LASSO_pos_SinglePenalty}
	\min_{\Z \in \R_+^{N \times L}} \| \A\Z - \Y \|_F^2 + \lambda \Rc_{\C} (\Z).
\end{align}


\section{Optimization via Gradient Descent} \label{sec:Optimization}
In order to approximate solutions to \eqref{eq:LASSO_pos_SinglePenalty} we use gradient descent. For $F(\Z) := \| \A\Z - \Y \|_F^2 + \lambda \Rc_{\C} (\Z)$ the gradient is given by
\begin{align*}
    \nabla F(\Z) = 2\A^T (\A\Z - \Y) + 2\lambda \; \vec^{-1} [\C \vv{\Z}],
\end{align*}
where $\vec^{-1}[\cdot] \colon \R^{NL} \rightarrow \R^{N\times L}$ inverts vectorization. To prevent gradient descent from leaving $\R_+^{N\times L}$, we replace the gradient by a projected version
\begin{align} \label{eq:ProjectedGrad}
    [\tilde{\nabla} F(\Z)]_{i,j} 
    := \begin{cases}
        0 & \text{ if } [\nabla F(\Z)]_{i,j} > 0,\ Z_{i,j} = 0, \\
        [\nabla F(\Z)]_{i,j} & \text{else.}
    \end{cases}
\end{align}
Note that $\tilde{\nabla} F(\Z)$ still points into a descent direction.
To compute a suitable step-size, we use the particular geometry of $\Rc_\C$. Let us define the descent ray function
\begin{align*}
    f_\Z (t) = F(\Z - t \tilde{\nabla} F(\Z)), \quad f_\Z \colon \R \rightarrow \R,
\end{align*}
for all $\Z \in \R^{N\times L}$. Note that $f_\Z$ can be written as
\begin{align*}
    f_\Z (t) 
    &= \left( ( \norm{\A \tilde{\nabla} F(\Z)}_F^2 +\; \lambda \vec^{-1} \left[\vv{\tilde{\nabla} F(\Z)}^T \C \vv{\tilde{\nabla} F(\Z)} \right] \right) t^2 \\
    &+ \left(2\inner{\A \tilde{\nabla} F(\Z), \A\Z - \Y}\phantom{(\tilde{\nabla} F(\Z) )^T} \right.\\
    &\phantom{+ (}\left.- \lambda \vv{\Z}^T \C \vv{\tilde{\nabla} F(\Z)} - \lambda \vv{\tilde{\nabla} F(\Z) }^T \C \vv{\Z} \right) t + c \\
    &= a t^2 + b t + c,
\end{align*}
where $c \in \R$ collects all terms not depending on $t$. We can now compute
\begin{align} \label{eq:sigma1}
    \tilde{\sigma}_1 = \inf_{\substack{Z_{i,j} > 0,\\ [\tilde{\nabla} F(\Z)]_{i,j} > 0}} \frac{Z_{i,j}}{\abs{[\tilde{\nabla} F(\Z)]_{i,j}}}
\end{align}
as the maximal allowed step-size to stay within $\R_+^{N\times L}$. Moreover, if $a > 0$, the function $f_\Z$ is strictly convex and the optimal (unconstrained) step-size is given by 
\begin{align} \label{eq:sigma2}
    \tilde{\sigma}_2 = -\frac{b}{2a}.
\end{align}
If $a \le 0$, we set $\tilde{\sigma}_2 = 1$. We thus use
\begin{align} \label{eq:sigmaClipping}
    \sigma = \min \{ \tilde{\sigma}_1, \tilde{\sigma}_2, 1 \}
\end{align}
as step-size for our algorithm. Note that the choice of $1$ as an upper bound for $\tilde{\sigma}_2$ and $\tilde \sigma_2$ is generic. Alternative choices would be $\frac{1}{\norm{\tilde{\nabla}F(\Z)}}$ or $\infty$.

\begin{figure}
\begin{algorithm}[H]
\caption{\textbf{:}  \textbf{Structural Sparse Recovery (SSR)}}
	\label{alg:GD}
	\begin{algorithmic}[1]
		\Require{$F(\Z) = \| \A\Z - \Y \|_F^2 + \lambda \Rc_{\C} (\Z)$ }
		\Let{$\X_0$}{$\0$}
		\While{stop condition is not satisfied}
		\Let{$\sigma_k$}{$\max \{ \tilde{\sigma}_1, \tilde{\sigma}_2, 1 \}$}
		\Comment{see \eqref{eq:sigma1}\ -\ \eqref{eq:sigmaClipping}}
		\Let{$\X_{k+1}$}{$\X_k + \sigma_k \tilde{\nabla} F(\X_k)$}
		\Comment{see \eqref{eq:ProjectedGrad}}
		\EndWhile
		\Statex
		\Return{$\X_\text{rec}$}
	\end{algorithmic}
\end{algorithm}
\end{figure}

\begin{remark}
    The structure encoding matrix $\C$ is a high-dimensional object. For an efficient implementation of Algorithm \ref{alg:GD} it is crucial to construct $\C$ such that it allows fast vector-matrix multiplication. The heuristic we propose in Appendix \ref{sec:C} fulfills this requirement, see Remark \ref{rem:Kronecker}.\\
    The recent work \cite{rick2017one} (which concentrates on $\X$ being a natural image) alternatively suggests to learn the projection operator onto the set of natural images via deep learning techniques and then to apply ADMM \cite{chambolle2011first} to reconstruct $\X$ from linear measurements. In our setting this would correspond to learning the projection onto $\mathbb S_\C^s$ for a certain type of ground-truths, e.g., seismic data, and then applying ADMM. Note, however, that compared to ours this approach has two major drawbacks. First, training a deep network requires massive amounts of data that are not always available. Second, one cannot interpret the structure of the sparsity in $\X$ from a learned network whereas this is possible to some extent when using $\C$, cf. Figure \ref{fig:ndt_alphabeta}. It would be interesting to compare both approaches numerically in future work.
\end{remark}

\subsection{Convergence of gradient descent} 

As the objective functional is non-convex, the question arises under which assumptions Algorithm \ref{alg:GD} can be expected to converge. The following observation provides a sufficient condition for this to happen. We omit the proof which is straight-forward (recall that $\C$ has only non-negative entries).

\begin{lemma}
    The objective function in \eqref{eq:LASSO_pos} is bounded from below on $\R_+^{N\times L}$. If
    \begin{align} \label{eq:Condition}
        \{ \Z \in \R_+^{N\times L} \colon \Rc_\C (\Z) = 0 \} \cap \ker(\A) = \curly{\0},
    \end{align}
    then it is coercive as well. Consequently, minimizers exist and all minimizers are contained in a finite ball around the origin.
\end{lemma}

\begin{remark}{}
    The assumption in \eqref{eq:Condition} is natural since it requires $\A$ to distinguish elementary structures from $\0$. When considering elementary structures generalizing row-sparsity, i.e., each column of the structured matrix contains at most one non-zero entry, the condition is fulfilled if $\A$ contains no zero columns. In case of more complex elementary structures with up to $s$-elements per columns in the structured matrix, \eqref{eq:Condition} is implied by $\A$ having a \emph{null-space property} of order $s$. In the context of solving underdetermined linear systems, null-space properties are a well-known concept from the compressed sensing literature \cite{foucart2013mathematical}.
\end{remark}


\section{Numerical Experiments}
\label{sec:Numerics}
Finally, let us demonstrate the power of the proposed Structural Sparse Recovery (SSR), see Algorithm \ref{alg:GD}. In Sections \ref{sec:ParameterTest} and \ref{sec:Comparison}, we perform basic tests on artificial data to validate our theoretical considerations. In Section \ref{sec:numerics_appl}, we present how SSR performs when applied to real-world data.

To quantify noise robustness, we use the notation of peak signal to noise ratio (PSNR). The reconstructed signal is denoted by $\X_{rec}$.

\subsection{Parameter heuristic test}
\label{sec:ParameterTest}

\begin{figure}
    \centering
    \subfloat[]{\includegraphics{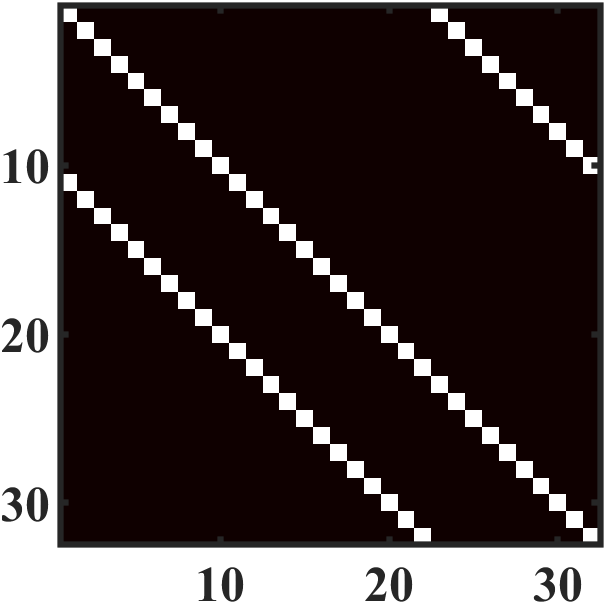}}\hfil
    \subfloat[]{\includegraphics{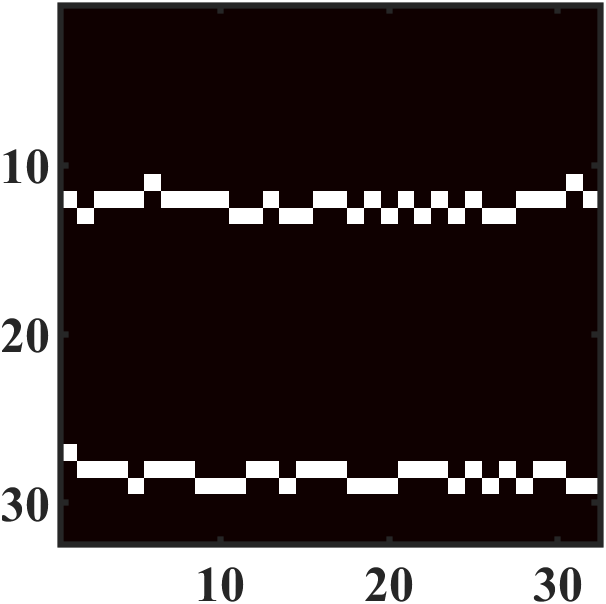}}\hfil
    \subfloat[]{
    \includegraphics{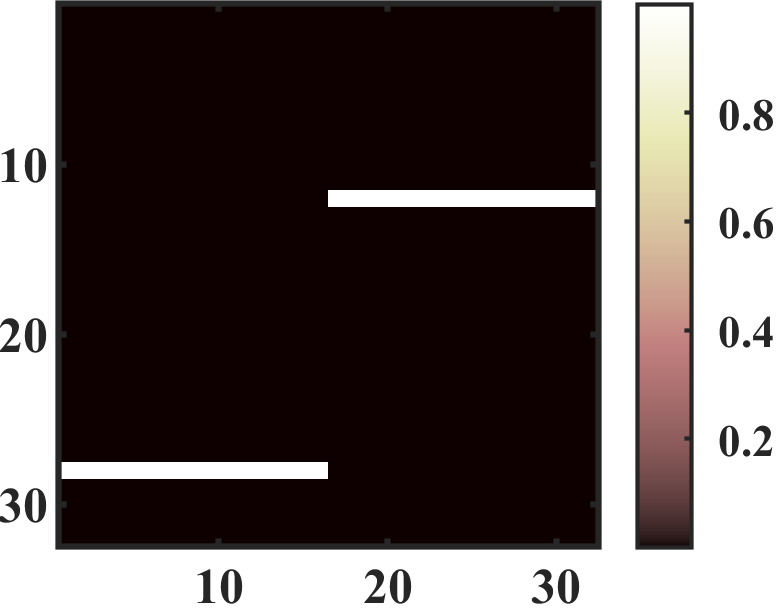}}
    \caption{Three examples of structures that may appear in applications: diagonals (a), oscillating lines (b), and partial row-sparsity (c).}
    \label{fig:X}
\end{figure}

\begin{figure}
    \centering
    \subfloat[]{
    \includegraphics{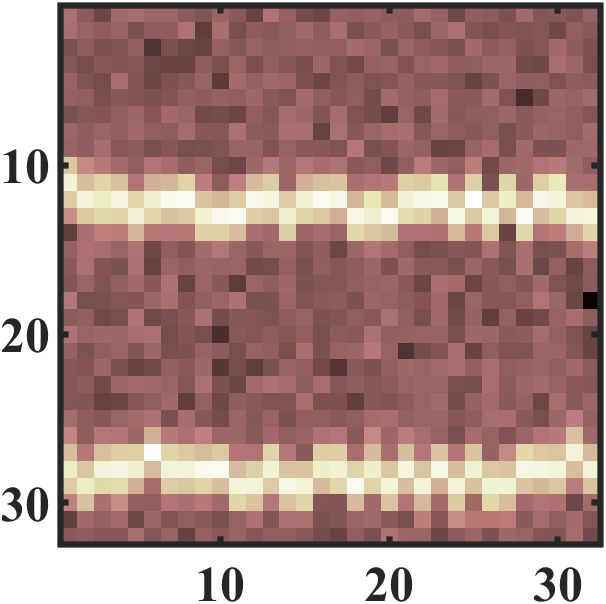}}\hfil
    \subfloat[]{
    \includegraphics{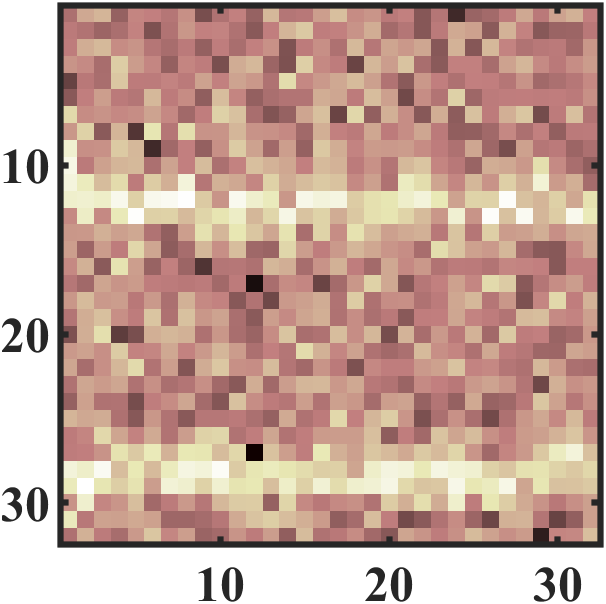}}\hfil
    \subfloat[]{
    \includegraphics{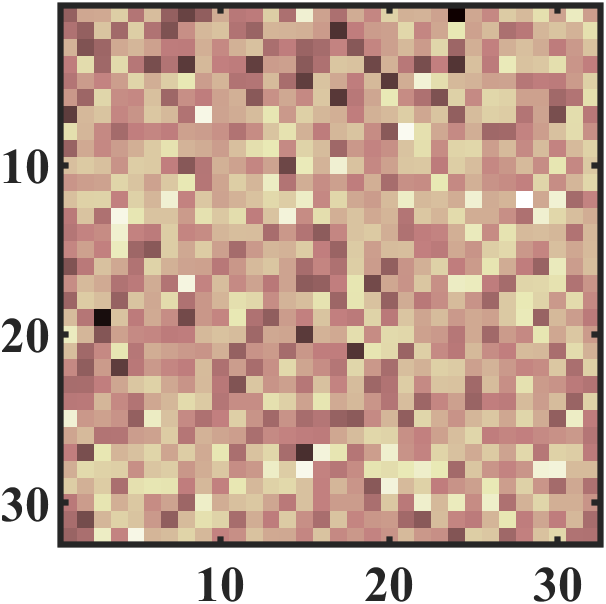}}
    \caption{Simulated data (oscillating lines) with low noise (a) (PSNR 24.92), medium noise (b) (PSNR 12.46), and high noise (c) (PSNR 1.97).}
    \label{fig:Measurements}
\end{figure}

Let us first empirically validate whether the heuristic for constructing $\C$ as proposed in Appendix \ref{sec:C} yields practical results and whether the computed estimates for the therein used parameters $\alpha$ and $\beta$ are meaningful. We choose three different $\X \in \mathbb{R}^{32\times32}$, see Fig. \ref{fig:X}, each of which contains two elementary structures (diagonal, slightly oscillating rows, partial rows), and a convolutional kernel $\A \in \mathbb{R}^{32\times32}$ modeling the measurement process. We add three different levels of Gaussian noise -- PSNR of $24.92$ (low noise), $12.46$ (moderate noise), and $1.97$ (high noise) -- to obtain corresponding measurements $\Y \in \mathbb{R}^{32\times32}$, see Fig. \ref{fig:Measurements}. When constructing $\C$, for fixed $\A$ and $\Y$ there are only finitely many thresholds $(\alpha,\beta)$ of interest, cf.\ \eqref{eq:Cbeta}. We can thus reconstruct $\X$ from $\A$ and $\Y$ for all possible choices of $\alpha$ and $\beta$ to obtain a benchmark for the heuristic. Note that we optimize in each reconstruction the additional parameters $\lambda_1,\lambda_2 \in [0,1]$, cf.\ \eqref{eq:LASSO_pos}, over a fine grid.

\begin{figure}
    \centering
    \includegraphics{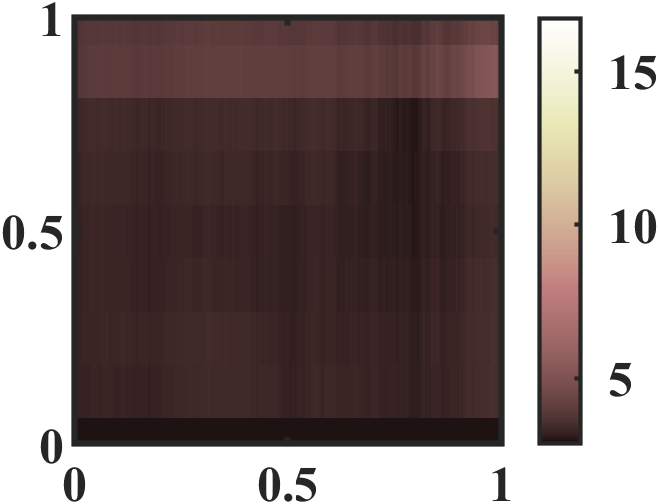}\hfil
    \includegraphics{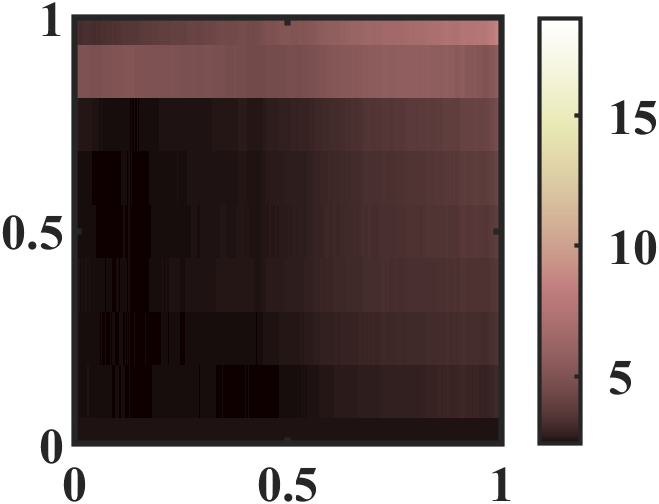}\hfil
    \includegraphics{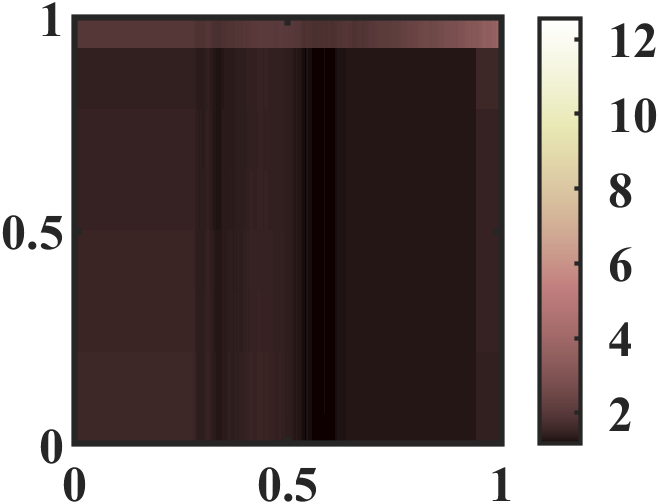}\\[0.35cm]
    \includegraphics{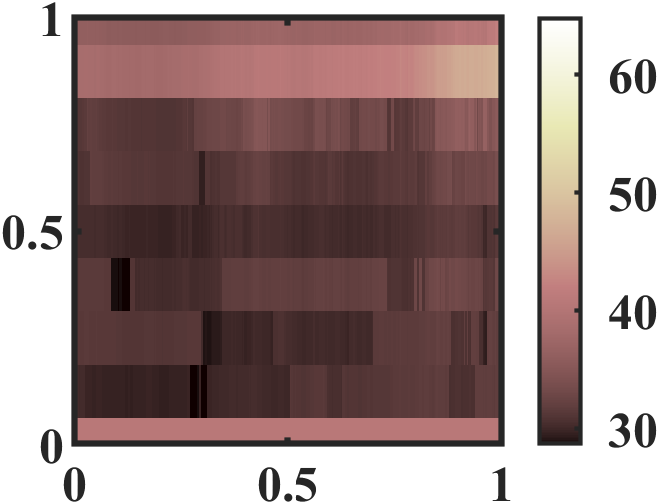}\hfil
    \includegraphics{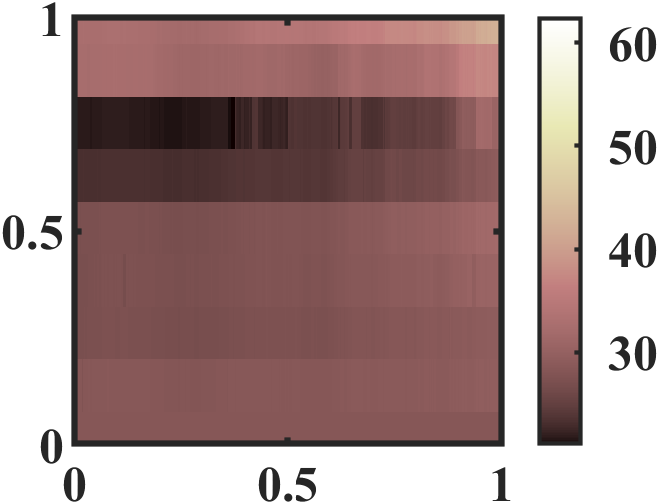}\hfil
    \includegraphics{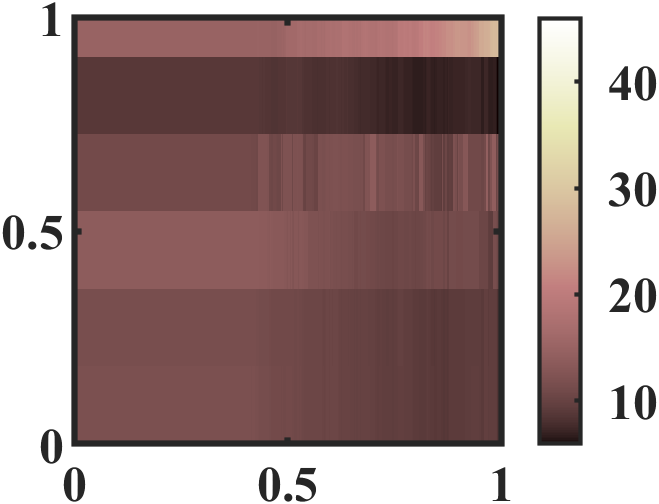}\\
    \subfloat[diagonal]{
    \includegraphics{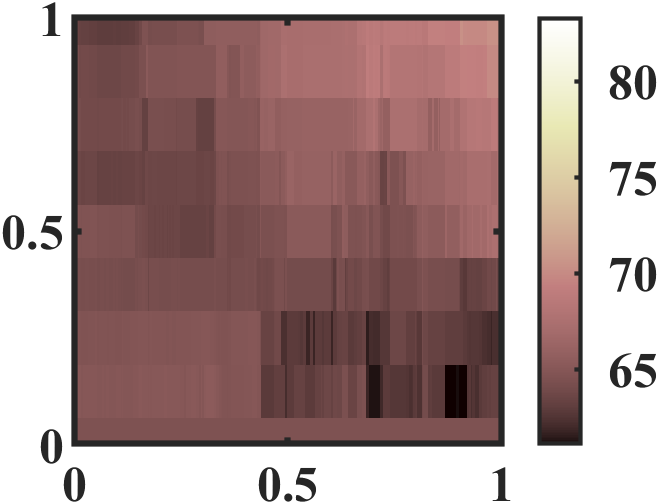}}\hfil
    \subfloat[oscillating]{
    \includegraphics{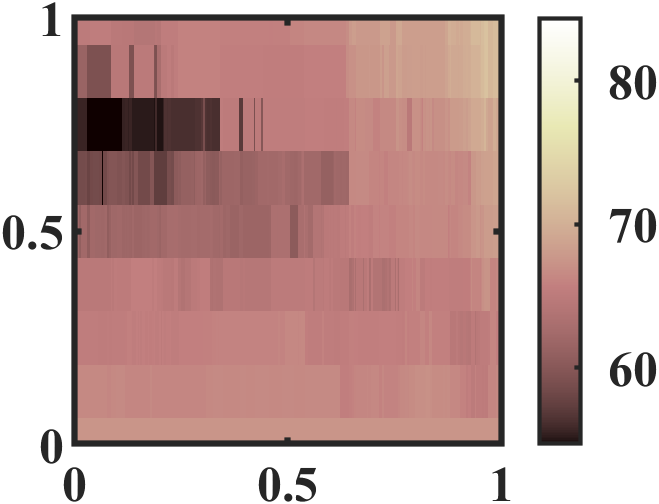}}\hfil
    \subfloat[row-sparse]{
    \includegraphics{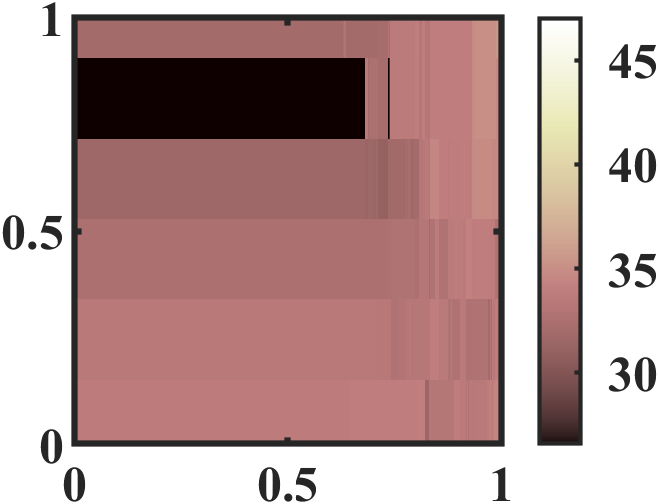}}
    \caption{Reconstruction error for different parameters $\alpha$ (y-axis), $\beta$ (x-axis) and different noise levels low (top row), medium (middle row), and high (bottom row).}
    \label{fig:ABgrid}
\end{figure}

Fig. \ref{fig:ABgrid} depicts the $\ell_1$-error after scaling $\X_\text{rec}$ to have the same $\ell_1$-norm as $\X$. This re-scaling balances the norm shrinkage caused by the penalty-term to allow for more meaningful comparisons. For the particularly chosen matrix $\A$, a choice of $\alpha=0.7313$ considers two atoms that are shifted by at most one pixel as similar while for $\alpha=0.99$ each atom is only similar to itself, i.e., we enforce row-sparsity. These thresholds can clearly be seen in Fig. \ref{fig:ABgrid} (in particular, when comparing a) and b) with c) in the low- and medium-noise cases). With increasing noise level it becomes more important to choose $\alpha$ sufficiently large in order to enforce structural sparsity. As the heuristic in Appendix \ref{sec:C} suggests, a choice of $\beta\leq\alpha$ yields optimal reconstruction results.
As the noise level increases, the value of $\beta$ has to be diminished such that more measurements are considered to be similar. This is as expected since more information is needed to reconstruct the signal in this case. Note that, for small to medium noise, the reconstruction performance is quite stable with respect to perturbations of the parameters $\alpha$ and $\beta$.

The diagonal case is most problematic in this regard (compare the color bars for the high-noise case) and demonstrates the limits of multiple measurements under this structure model. The strength of MMV is, to consider information of several columns at once to reconstruct the support. In the extreme case of row-sparsity, all columns can be taken into account at once. For localized structures such as Fig. \ref{fig:X} a) where the support changes fast the advantage of MMV is limited as only a small neighborhood carries information about the actual support of a column. Considering a larger neighborhood (i.e., decreasing $\beta$ in our model) can help to improve the results. However, due to the fast changes of the support, we also need to consider more atoms as similar in the same step (i.e., decreasing $\alpha$) which diminishes the gained information. Thus, if the noise level is too high, the model is not able to gain enough information for a suitable reconstruction. One way to overcome this problem might be to design a matrix $\C$ that is highly adapted to allow only diagonal structures.

\subsection{Comparison}
\label{sec:Comparison}

Let us now benchmark SSR against five state-of-the-art methods for the detection of sparsity and row-sparsity: Iterative Soft-Thresholding Algorithm (ISTA) \cite{daubechies2004iterative}, Group-Iterative Soft-Thresholding Algorithm (G-ISTA) \cite{yuan2006model}, Orthogonal Matching Pursuit (OMP) and Rank Aware Order Recursive Matching Pursuit (RA-ORMP) \cite{Tropp06a,davies12}, and Sparse Bayesian Learning for row-sparsity (SBL) \cite{Wipf07} (where we use the implementation of Z. Zhang \cite{sblWebpage}). ISTA is a technique for sparse recovery and G-ISTA its generalization to group sparse recovery of signals from underdetermined linear systems. 
RA-ORMP is an extension of the OMP algorithm for simultaneous sparse approximation. Last but not least, Sparse Bayesian Learning is a probabilistic regression method that has been developed in the context of machine learning.

We use for $\X$ the same structure types as considered in the previous section, see Fig. \ref{fig:X}. Note that the different methods intend to minimize different sparsity measures. While OMP and ISTA minimize the individual sparsity of each column, G-ISTA, RA-ORMP and SBL minimize the row-sparsity. Hence, for different $\X$ some algorithms might be more suitable than others. Table \ref{tab:sparsity} shows that Fig. \ref{fig:X} a) is clearly not row-sparse and thus favors our approach as well as the single measurement algorithms. In Fig. \ref{fig:X} b) row sparsity and structural sparsity are of comparable order. Finally, Fig. \ref{fig:X} c) has the same row and structural sparsity such that the comparison is not biased towards one single MMV method. When applying the different reconstruction algorithms we assume that the sparsity level is known, i.e., the iterative algorithms perform the exact amount of iterations needed. Moreover, the exact rank of $X$ was given to RA-ORMP. Further hyper-parameters have been optimized over a grid to ensure best possible performance of all methods in the direct comparison.

\begin{table}[]
    \centering
    \begin{tabular}{|c|c|c|c|}
    \hline \\
    Structure & max column & row-sparsity & structural \\
    Fig. \ref{fig:X} & sparsity $\max\limits_k\|\x_k\|_0$ & $\|\X\|_{\text{row}-0}$ & sparsity $\|\X\|_{\C,0}$ \\
    \hline \\
    a) & 2 & 32 & 2\\
    b) & 2 & 6 & 2\\
    c) & 1 & 2 & 2\\
    \hline
    \end{tabular}
    \caption{Sparsity of the structures shown in Fig. \ref{fig:X} under different norms.}
    \label{tab:sparsity}
    \vspace{-20pt}
\end{table}

\begin{figure}
    \centering
    \subfloat[diagonals]{
    \includegraphics{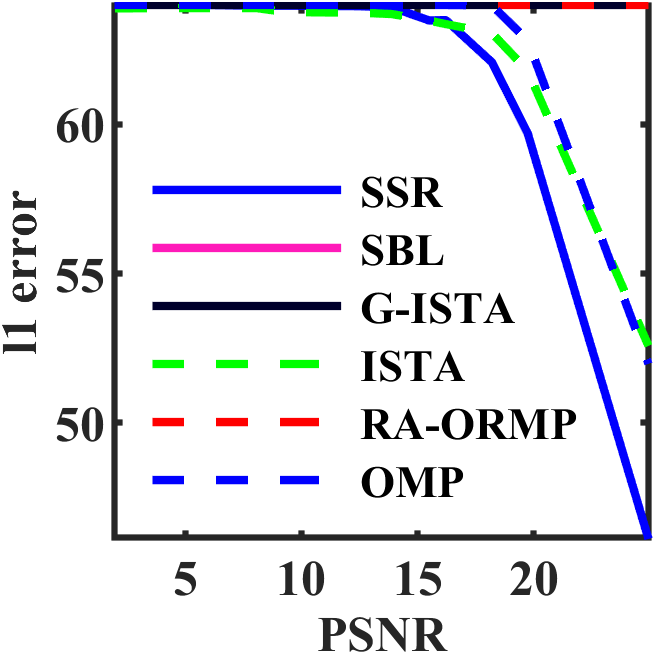}}\hfil
    \subfloat[oscillating lines]{
    \includegraphics{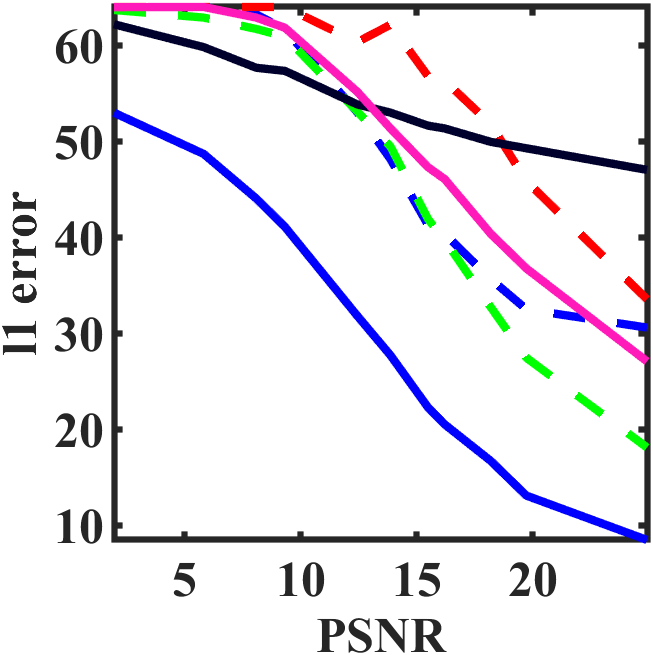}}\hfil
    \subfloat[row-sparse]{
    \includegraphics{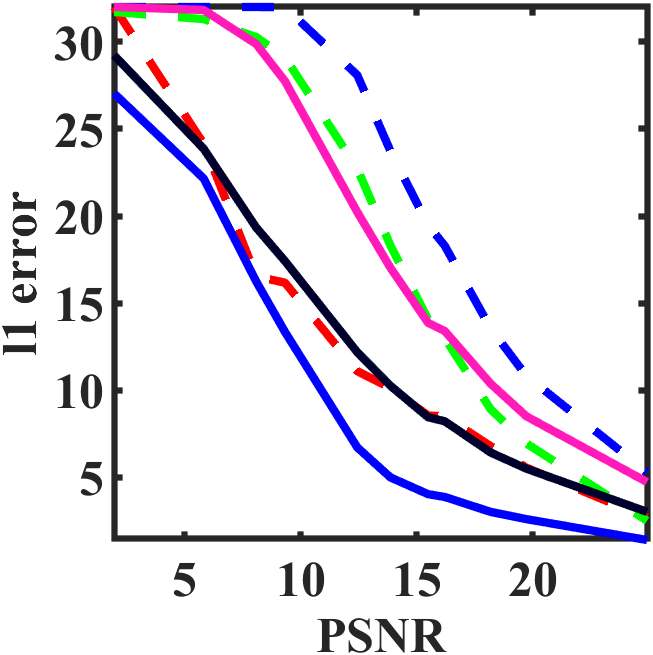}}
    \caption{Mean $\ell_1$-error (scaled) vs. different noise levels for the different structures.}
    \label{fig:benchmark}
\end{figure}

\begin{figure}
    \centering
    \includegraphics{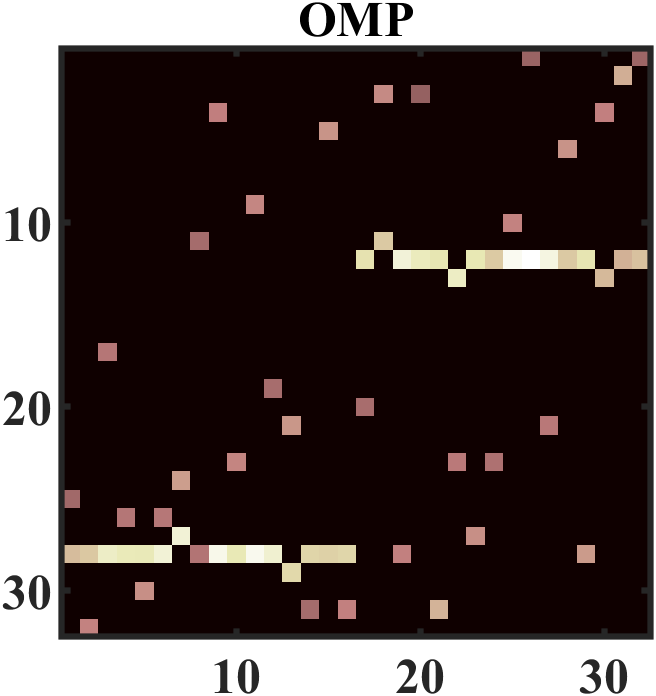}
    \includegraphics{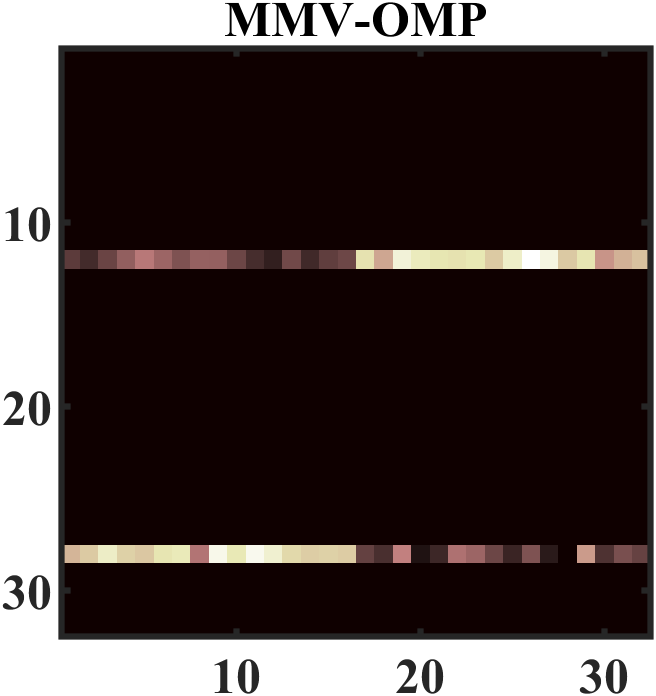}
    \includegraphics{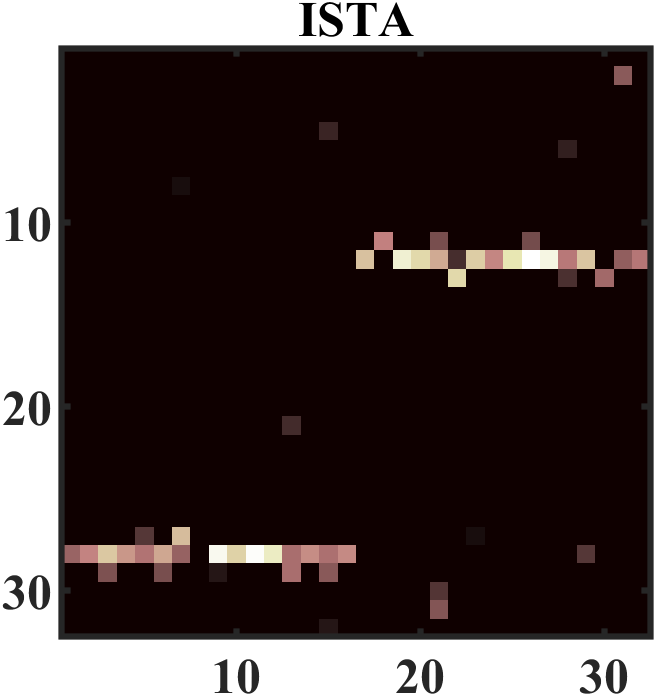}\\[0.25cm]
    \includegraphics{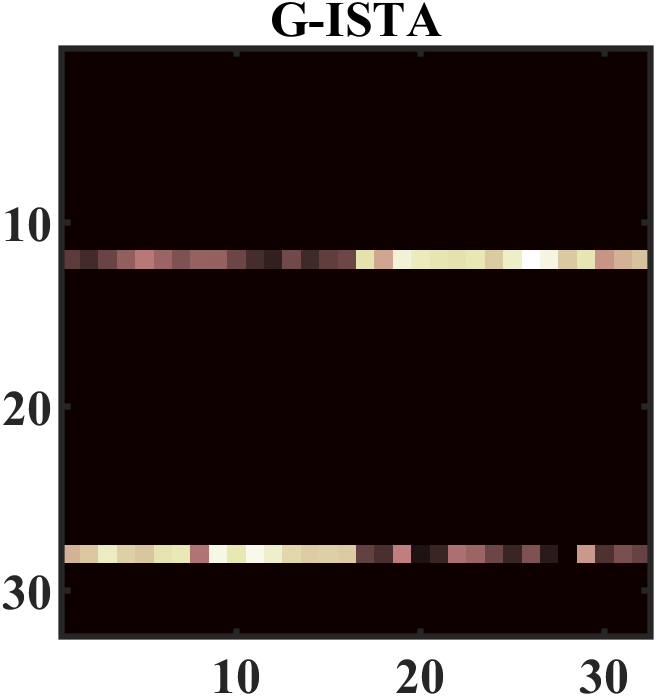}
    \includegraphics{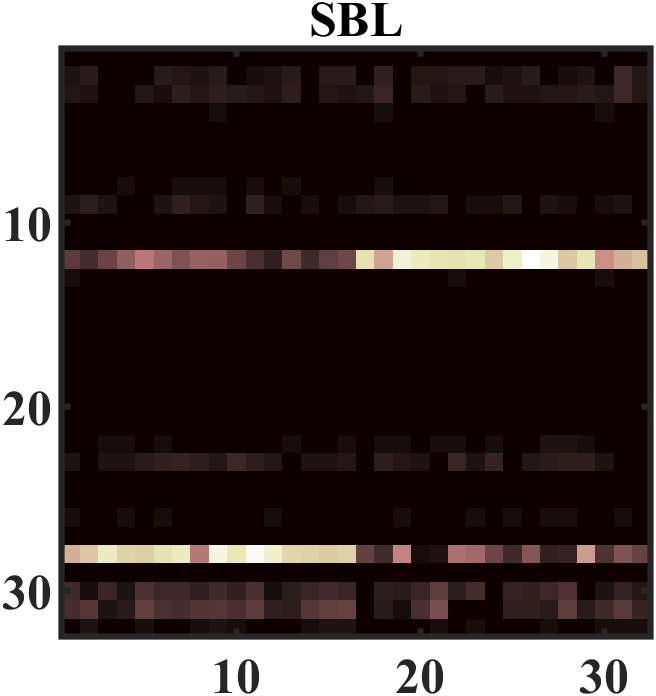}
    \includegraphics{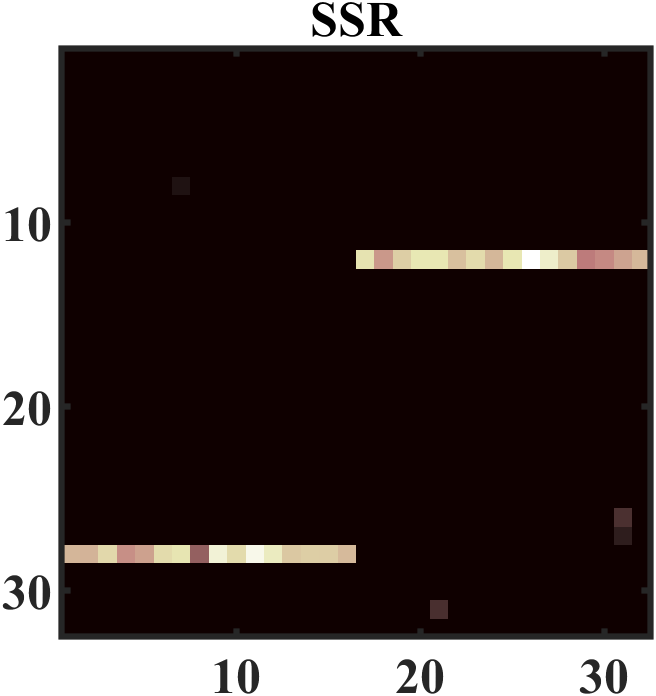}
    \caption{Absolute value of reconstructions for different algorithms on row-sparse structures with medium noise level (PSNR $12.46$).}
    \label{fig:recon_comparison}
\end{figure}

To create our test data we use four different measurement matrices $\A$, two convolution matrices $\A \in \R^{32\times32}$ with a Gauss kernel (similar to the previous section) as well as two random Gaussian matrices $\A \in \R^{10\times32}$. This is a typical case of undersampling as considered in compressed sensing. In order to compare the results we re-scale the reconstruction $\X_{rec}$ as described above. In Fig. \ref{fig:benchmark} the average $\ell_1$-error of $300$ runs is plotted for $11$ different noise levels (PSNR of 2 to 25) where the range of the y-axis is restricted by the $\ell_1$-norm of the ground-truth. It is clearly visible that the SSR algorithm outperforms all competitors, especially in the case of oscillating lines. This suggests that it is more robust to row-sparsity defects. A reconstruction example for each algorithm on the row-sparse structure is depicted in Fig. \ref{fig:recon_comparison}. We display the reconstructions in absolute value since small negative entries would otherwise change the colormap where $0$ is supposed to be black, and thus make a comparison with the original more difficult.

\subsection{Applications}
\label{sec:numerics_appl}

Finally, we apply SSR to real-world data to show its capability of analyzing mixtures of different types of complex structures.

\subsubsection{Non-destructive testing}

\begin{figure}
    \centering
    \subfloat[]{
    \includegraphics{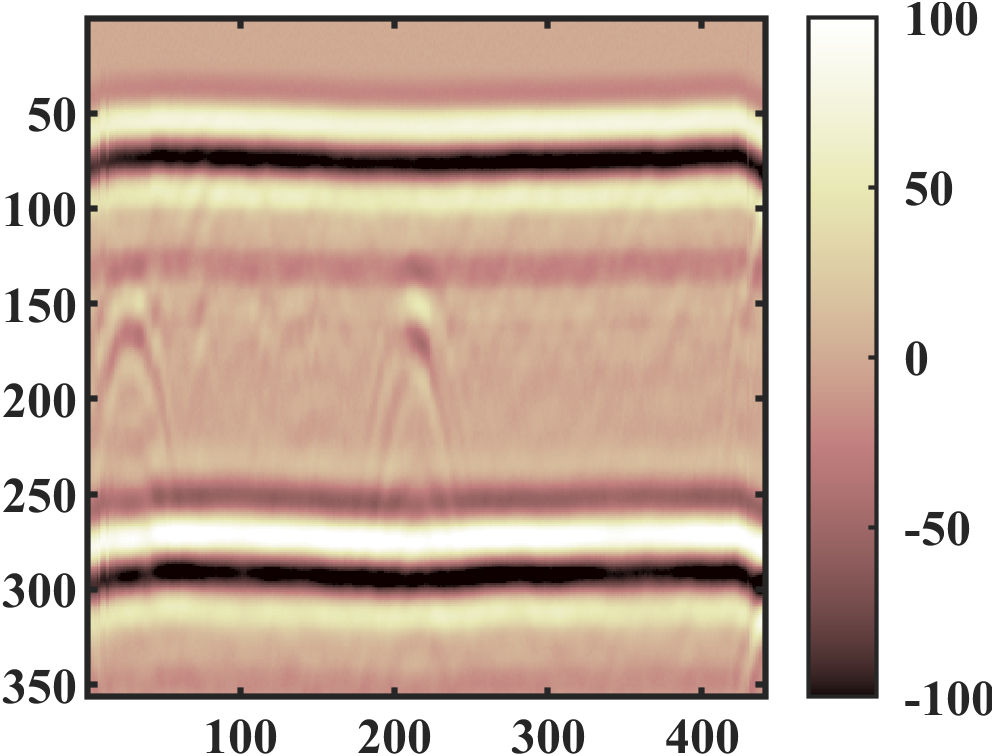}}
    \label{fig:ndt1}\hfil
    \subfloat[]{
    \includegraphics{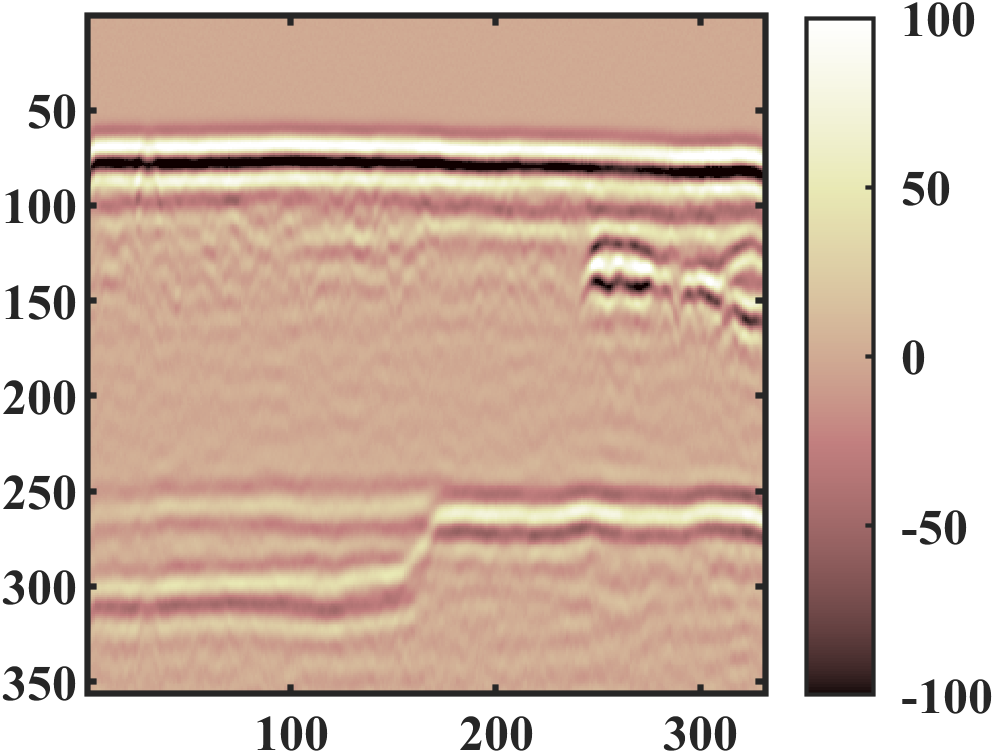}}
    \label{fig:ndt2}
    \caption{Original ultrasonic non-destructive testing measurements of pores (a) and lack of fusion (b).}
    \label{fig:ndt_org}
\end{figure}

\begin{figure}
    \centering
    \subfloat[]{
    \includegraphics[align=c]{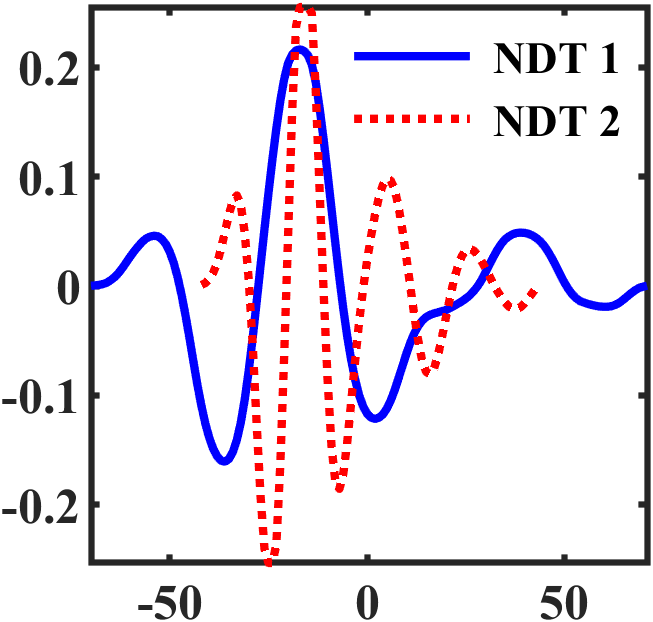}}
    \label{fig:ndt_impulse}\hfil
    \subfloat[]{
    \includegraphics[align=c]{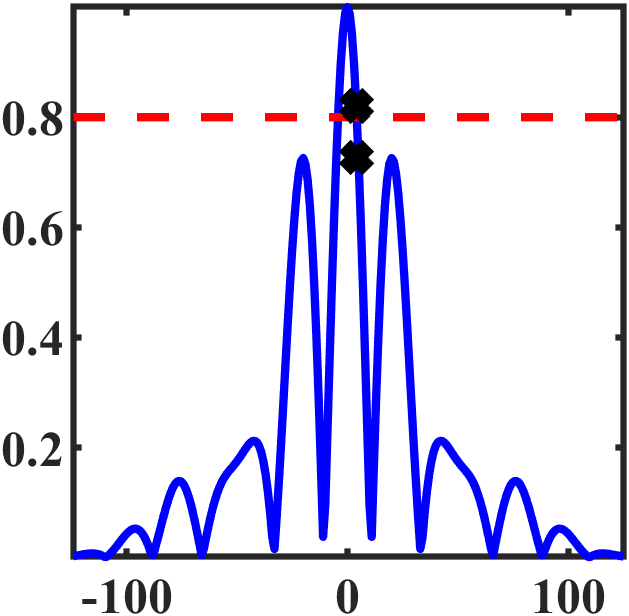}}
    \label{fig:ndt1_alpha}\hfil
    \subfloat[]{
    \includegraphics[align=c]{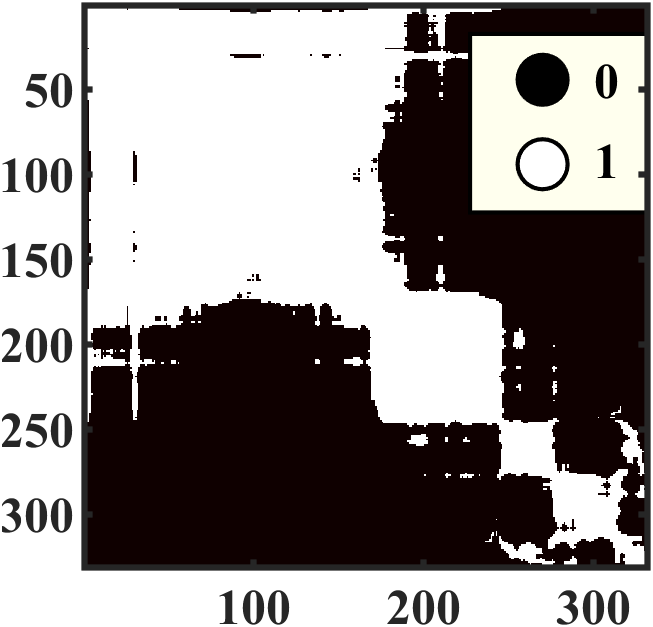}}
    \label{fig:ndt2_beta}
    \caption{Impulse functions (a), auto-correlation and the parameter $\alpha$ (dashed line, separates the fourth and fifth entry marked in black) (b) exemplary for NDT 1 in (a), and matrix $\C_{\beta}$ ($\beta=0.75$) for the data in Fig. \ref{fig:ndt2}.}
    \label{fig:ndt_alphabeta}
\end{figure}

The first example comes from the manufacturing industry in the field of non-destructive testing. Here, one tries to detect material defects and other anomalies from ultrasonic images of an object. In Fig. \ref{fig:ndt_org} ultrasonic images from scanning the weld seam of a steel pipe are depicted where the different structures result from the lateral signal and the back wall echo (horizontal lines) and the anomalies. The two images show different types of anomalies in the weld: in Fig. \ref{fig:ndt_org} (a) we presumably see pores, whereas Fig. \ref{fig:ndt_org} (b) shows a lack of fusion at the end of the pipe, where the last part of the weld seam has been ground. For representation of a received signal, one supposes that it can be obtained as a linear combination of time-shifted, energy-attenuated versions of the reconstructed pulse function (see Fig. \ref{fig:ndt_alphabeta} (a)), where each shift is caused by an isolated flaw scattering the transmitted pulse \cite{Bossmann12}. This linear combination is described by the measurement matrix $\A$, a convolution matrix. 

We expect a sparse reconstruction $\X_{rec}$ as we assume the material to have only few anomalies. Moreover, $\X_{rec}$ exhibits structure since adjacent measurements will give similar results. In particular, any anomaly will be seen in several adjacent measurements.

\begin{figure}
    \centering
    \subfloat[]{
    \includegraphics{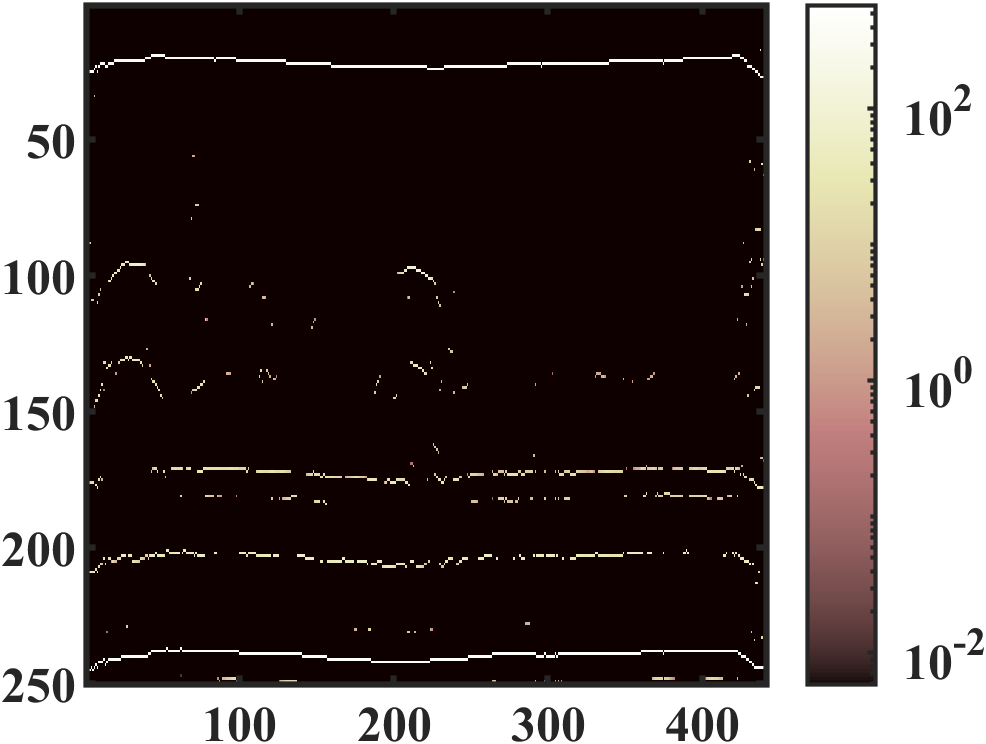}}
    \label{fig:ndt1_rec}\hfil
    \subfloat[]{
    \includegraphics{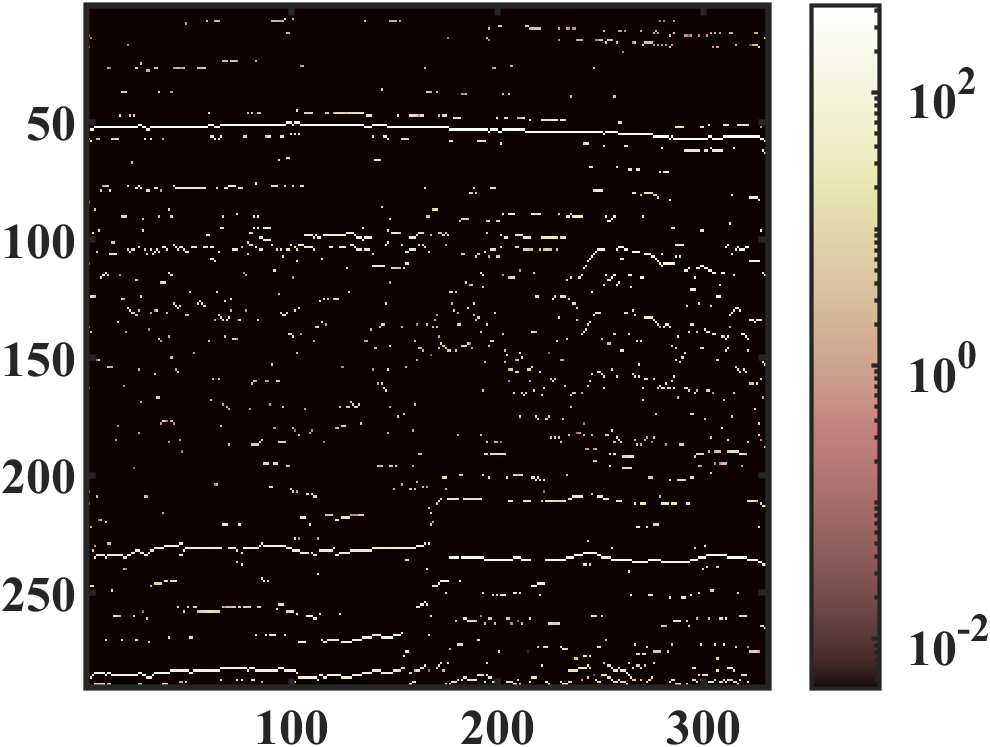}}
    \label{fig:ndt2_rec}
    \caption{Reconstructed structures (absolute value) in log-scale.}
    \label{fig:ndt_reconstruction}
\end{figure}

For the reconstruction of $\X$ we use the SSR model \eqref{eq:LASSO_pos} with $\lambda_1=10^{-4}$ and $\lambda_2=10^6$. (We choose here $\lambda_1 = \mathcal O (L^{-2})$ following the intuition in Remark \ref{rmk:parameterchoice}. Since the data is rather noisy we set $\lambda_2 = \mathcal{O}(L^2)$ instead of $\lambda_2 = \mathcal{O}(1)$ to prevent SSR from detecting multiple similar atoms per column.) The parameter $\alpha$ can be estimated using the auto-correlation of the wave impulse. Given the ultrasonic speed in the material, the time sampling of the signal, and the measurement setup, we can derive by how many pixels a signal from the same source can shift in between different measurements. For the examples given here a maximum shift of four respectively two pixels indicates a signal coming from the same source. Now, we choose $\alpha$ such that it separates the fourth and fifth respectively the second and third entry of the auto-correlation, cf. Fig. \ref{fig:ndt_alphabeta} (b).
For the example with the pores we set $\alpha=0.8$ and for the example with the lack of fusion $\alpha=0.78$. For the construction of $\C_2$ we choose a smaller value, namely $\alpha/5$, to enforce a gap between different structures. This way, the reconstruction results are further improved.
As the noise level is low, we set $\beta=0.75$ in both examples. We can see that for this choice the data structure already becomes evident in $\C_\beta$, e.g., the blocks appearing in Fig. \ref{fig:ndt_alphabeta} (c) divide the measurements (columns in Fig. \ref{fig:ndt_alphabeta} (c)) in three qualitatively different segments: the first segment contains two bands of which the lower one is diffuse, the second segment contains two well-separated bands, and the third contains the same two bands plus an additional artifact. The example thus perfectly fits into the structural sparsity framework of the paper.
The reconstruction results are shown in Fig. \ref{fig:ndt_reconstruction}. It can be seen that the method is able to detect both types of defects.

\subsubsection{Seismic exploration}

The second example is from the field of seismic exploration where the goal is to reconstruct soil layers in order to find natural resources such as gas and oil deposits or minerals. For this purpose, a wave is generated in a field experiment (e.g. by an explosion) and sensors that are arranged on a grid around the source register incoming seismic waves over time. The measured data form a 3D tensor of which we are looking at a slice (measured data along a straight line). For this reason, the example is very similar to the previous one and the reconstruction can be determined again with the same approach provided that the wave function is known (analogous to the pulse function in non-destructive testing). In \cite{Bossmann15,Bossmann16} it is described how to calculate a wave model as accurate as possible from the measurement data. This wave is then used again to set up the measurement matrix $\A$.  

\begin{figure*}
    \centering
    \subfloat[]{
    \includegraphics{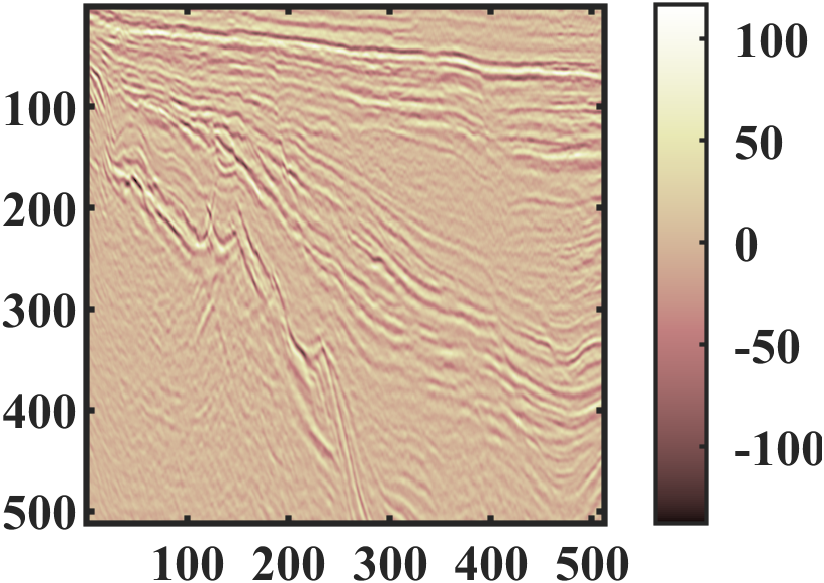}}
    \label{fig:seismic_org1}\hfil
    \subfloat[]{
    \includegraphics{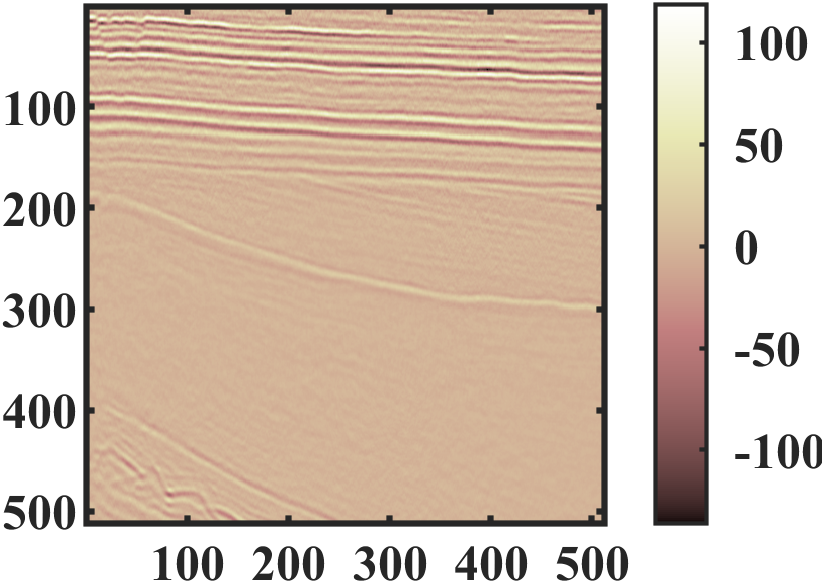}}
    \label{fig:seismic_org2}\hfil
    \subfloat[]{
    \includegraphics{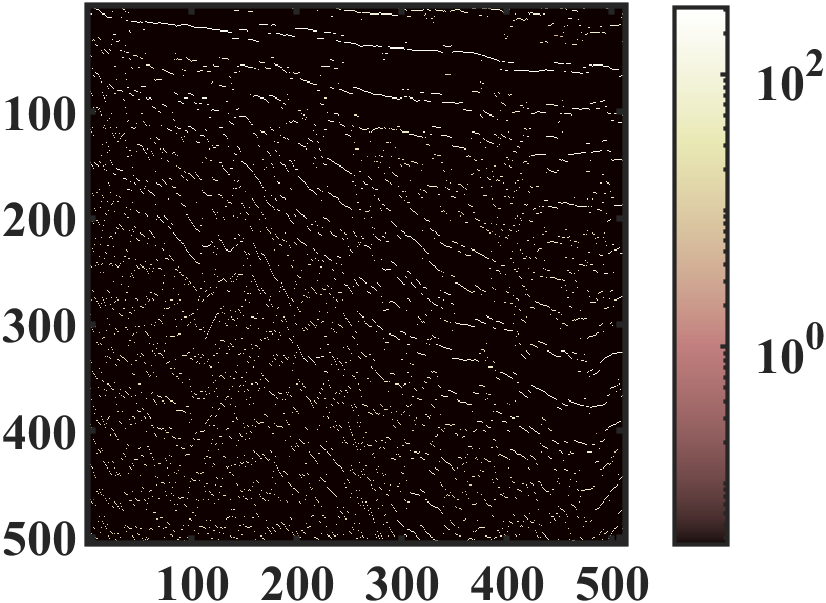}}
    \label{fig:seismic_recon1}\hfil
    \subfloat[long and short wave]{
    \includegraphics{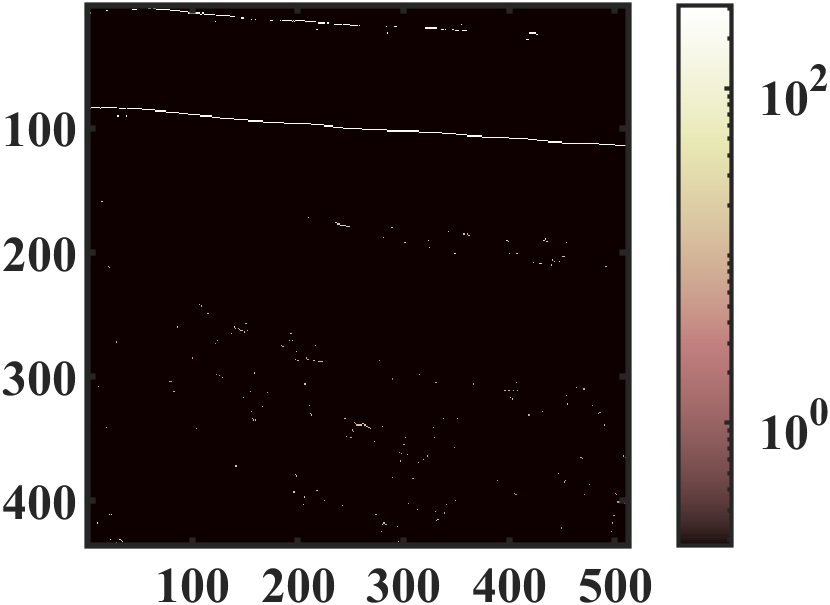}
    \includegraphics{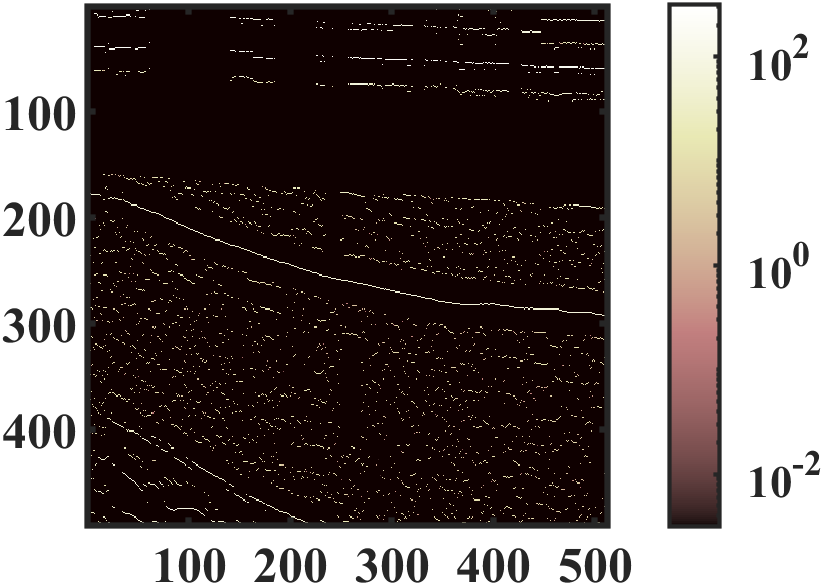}}
    \label{fig:seismic_recon2}
    \caption{Original seismic measurements (a,b) and the reconstructed structures (absolute value) in log-scale (c,d). Note that we show the reconstruction for the long and short wave impulse separately in (d).}
    \label{fig:seismic_org}
\end{figure*}

In Fig. \ref{fig:seismic_org} one can see slices from two different experiments with very different structures. While the seismic waves in the left image are all quite similar, the right image contains a mixture of longer and shorter waves.
This indicates that for the right image the combination of a shorter and a longer wave impulse may be useful. A short wave impulse would reconstruct the diagonal structures in the middle and bottom left of the image well, but would not detect all horizontal structures in the upper part. This could only be achieved by using many short waves contradicting the assumption of sparsity. If we used only one longer wave instead, we would not capture the diagonal structures.

The parameters can be estimated as in the previous section. To compensate noise and model approximation we choose $\alpha\cdot 10^{-3}$ when constructing $\C_2$, i.e., no interference between structures is admitted. The values of the parameters are depicted in Table \ref{tab:seismic}. When choosing $\lambda_1$ and $\lambda_2$ we apply in Example 1 the same reasoning as in the non-destructive testing experiment; in Example 2, we increase both parameters to prevent long seismic waves from being replaced by multiple short ones.

\begin{table}
  \centering
\begin{tabular}{lcccc}
            & $\alpha$ & $\beta$ & $\lambda_1$ & $\lambda_2$\\
\hline
Example 1   & 0.5    & 0.4   & $10^{-6}$   & $10^{6}$ \\
Example 2   & 0.7    & 0.6   & 1         & $1.5*10^{6}$\\
\end{tabular}
\caption{Parameters for seismic exploration.}
\label{tab:seismic}
    \vspace{-20pt}
\end{table}

The results are given in Fig. \ref{fig:seismic_org}. Again, all kinds of structures have been reconstructed. The usage of two wave functions pays off as the row-like structures (Fig.~\ref{fig:seismic_org}~(d), left) could be reconstructed separately from the diagonal-like structures. 

\subsubsection{Meteorology}\label{sec:meteor}

\begin{figure*}
    \centering
    \subfloat[]{
    \includegraphics[align=c]{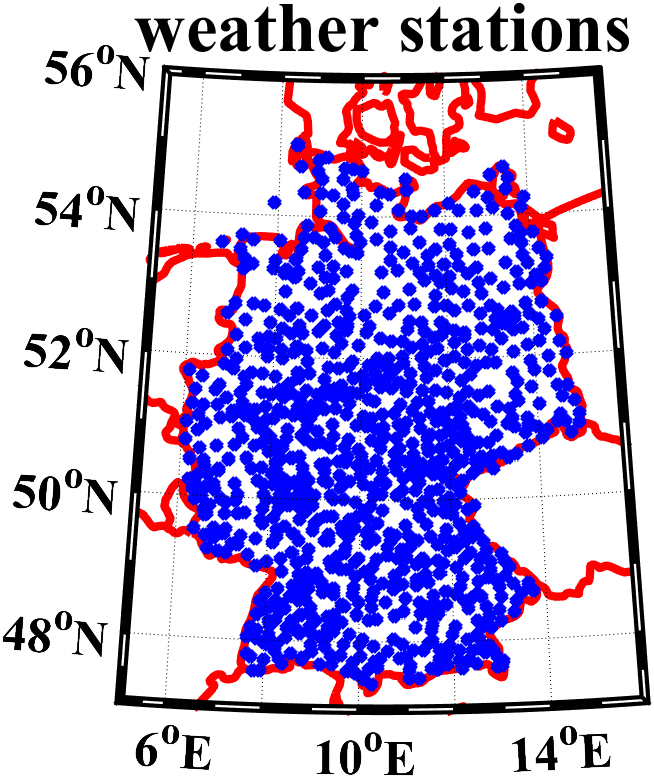}}
    \label{fig:meteo_stations}\hfil
    \subfloat[]{
    \includegraphics[align=c]{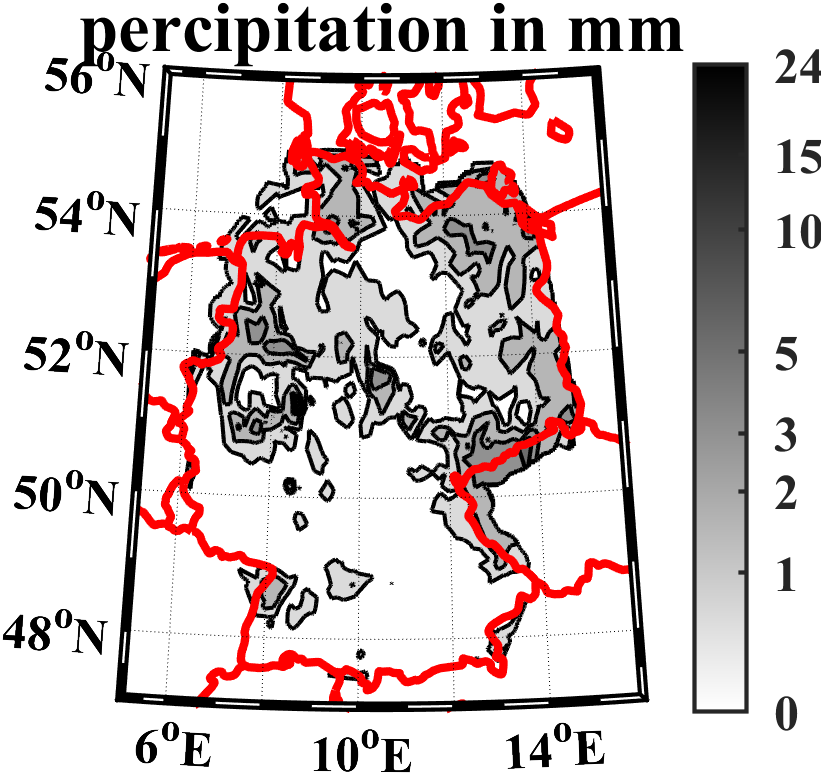}}
    \label{fig_meteo_rain}\hfil
    \subfloat[]{
    \includegraphics[align=c]{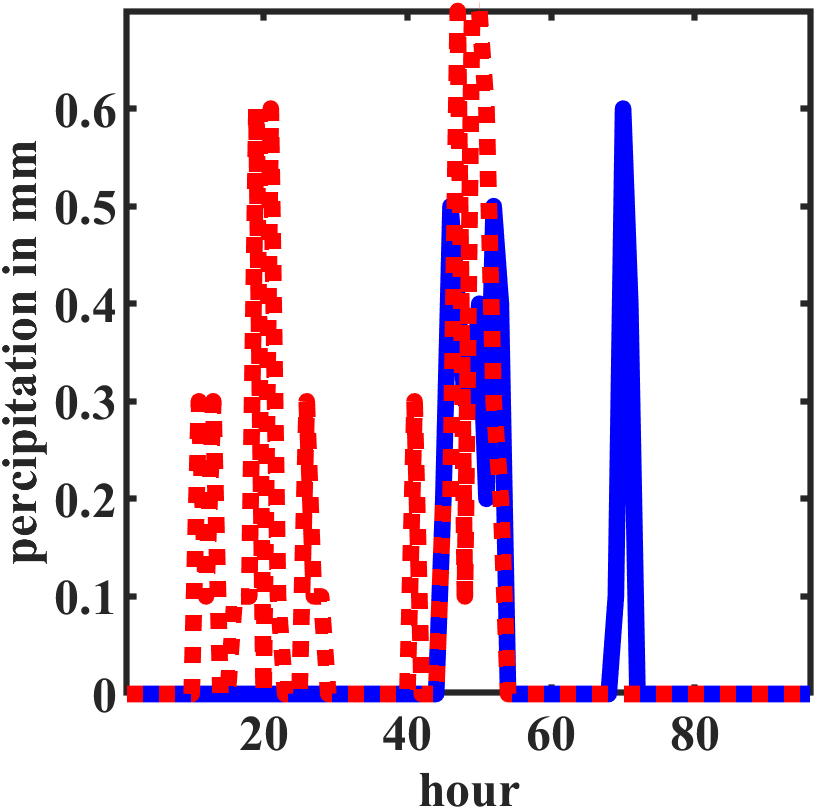}}
    \label{fig_meteo_staionary}\hfil
    \subfloat[]{
    \includegraphics[align=c]{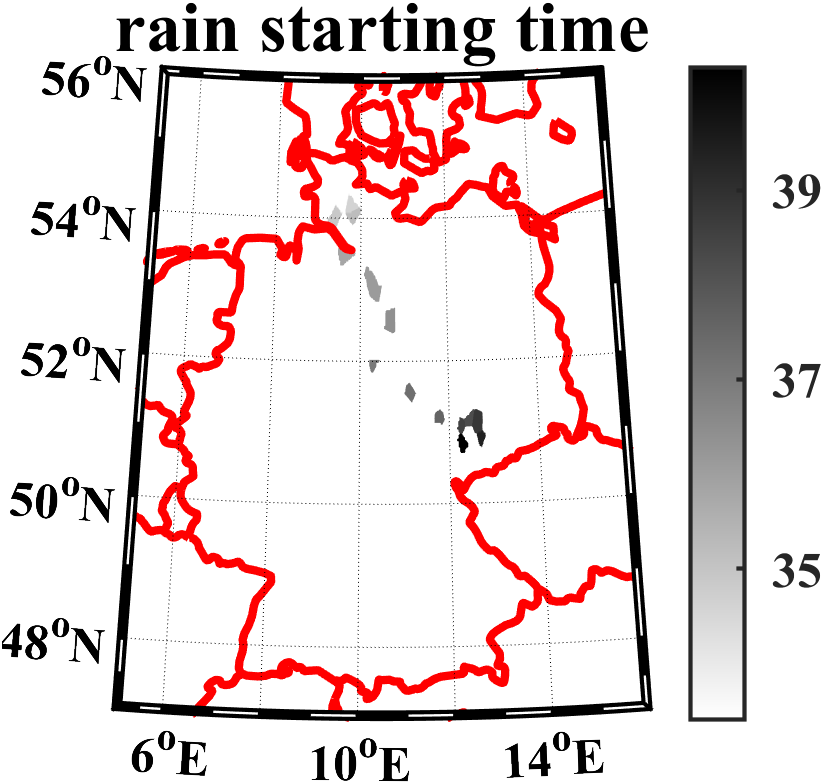}}
    \label{fig:weather_results1}\hfil
    \subfloat[]{
    \includegraphics[align=c]{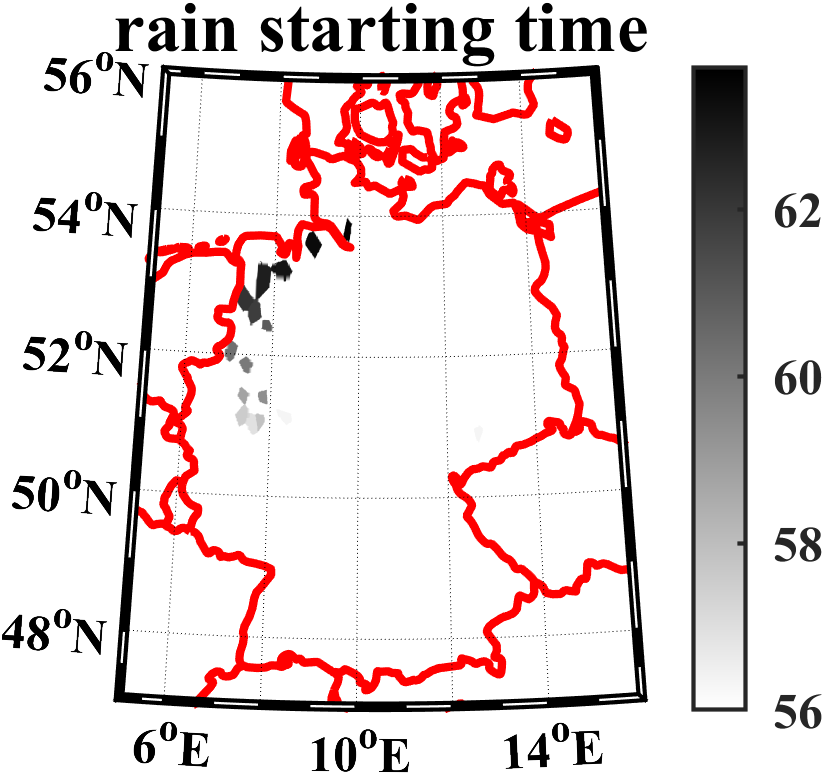}}
    \label{fig:weather_results2}
    \caption{Overall precipitation in Germany over 96 hours in November 2008 (b) (log scale) measured at locations (a). Exemplary data of two weather stations is given in (c). Tracked rain showers with longest duration (d,e).}
    \label{fig:meteo}
\end{figure*}

As a third variation, we analyze the hourly precipitation in Germany from 25th to 28th of November 2008 (96 hours) using data of 932 weather stations shown in Fig. \ref{fig:meteo} (a). Note that stations that were moved during this time or had too many missing values have not been taken into account. From these data we want to extract rain areas that are connected in time and space over as long a period as possible. We assume that the wind speed is constantly less than 75 km/h.
Translated into our setting, this means that a connected rain area is a structure and $\X \in \mathbb R^{465\times932}$ divides the total rain into disjoint, connected rain areas. The atoms of the system matrix $\A$ will represent rain showers of different lengths (two to six hours) observed at one station. Hence, $\A\in \{0,1\}^{96\times 465}$ is the block band matrix 
\begin{align*}
\A = [\A^2,\A^3,\ldots,\A^6] \text{ and }
(A^k)_{j,l} = \begin{cases}
1   & 1\leq j-l<k \\
0   & \text{else}
\end{cases}.
\end{align*}
Instead of correlation between the atoms we use the starting time difference of these showers and instead of correlation between measurements we use the distance between the weather stations. This way, $\C\in\{0,1\}^{932\cdot465\times932\cdot465}$ encodes the geographical information. More precisely, 
$C_{(j,l),(j',l')}=0$ only if the distance between the stations is smaller than $75$km/h times the starting time difference. Note that there is more meteorological data, such as wind direction or precise wind speeds of the according days available that could be used to construct a refined $\C$ such as wind direction or precise wind speeds of the according days. However, for this example we stick to the simpler model.

The measurement data shows the amount of precipitation in mm per hour. We simplify this to a binary matrix by setting $\Y\in\mathbb \{0,1\}^{96\times 932}$ where 0 indicates no rain and 1 indicates rain at a certain station at a certain time. This perfectly suits our system matrix $\A$ and guarantees that the linear system has a solution.

Despite the tremendous size of this system of equations it can be solved efficiently due to the fact that $\C$ is binary and can be factorized as a Kronecker product. Thus, the matrix multiplications can be reduced to summations. Moreover, we do not have to take the detour described in \eqref{eq:Extended} as all involved matrices are non negative by definition. 

In our computations we set $\lambda_1=200$ to just get one rain area per time interval and $\lambda_2=0$ since we do not enforce the structures to be sparse.

In Fig. \ref{fig:meteo} the two longest connected rain areas during the recorded period are depicted. In Fig. \ref{fig:meteo} (d) we see a rain area moving from north to south-east starting at about 34 hours and ending at 51 and in Fig. \ref{fig:meteo} (e) we see a rain area at the boarder to the Netherlands moving from south to north and starting at 56 hours and ending at about 75.


\section{Conclusion}

We presented in this paper a novel and user friendly model for handling structural sparsity in jointly sparse reconstruction tasks like MMV. We showed that our model is capable of representing various sparsity patterns of practical relevance. It, furthermore, lends itself to efficient reconstruction via projected gradient descent. In practice, the model parameters that determine the structures of interest can be heuristically learned from application specific environmental parameters, which are usually known to the practitioner. Comprehensive empirical studies confirmed efficacy of our method.\\
There remain several intriguing questions to be addressed in future work. First, we restricted the present work to introducing resp.\ motivating the model and empirically demonstrating its applicability. It would be desirable to derive supportive theoretical guarantees, both for the parameter heuristic in Appendix \ref{sec:C} and the (non-convex) projected gradient descent.\\
Second, it might be beneficial to modify the construction of $\C$ to allow entries with values in $(0,1)$, e.g., by using a soft-thresholding procedure. This could be used to enlargen the scope of representable structures, although it will presumably complicate parameter tuning. \\
Finally, it would be worthwhile to examine alternative constructions of $\C$ and compare their performance to our current heuristic in Appendix \ref{sec:C}.


\appendix

\subsection{Structural Sparsity Matrices} 
\label{sec:C}

There is no unique admissible way to construct a structural matrix $\C$ for \eqref{eq:Scs}. For our purpose a simple heuristic approach suffices that constructs $\C$ from $\A$ and $\Y$. As our numerical examples in Section \ref{sec:numerics_appl} illustrate, the approach lends itself to further refinement depending on the concrete application. It is based on the idea that if two measurement vectors $\y_l = \A \x_l$ and $\y_{l'} = \A \x_{l'}$ are similar, then the corresponding atoms activated by $\x_l$ and $\x_{l'}$ in $\A$, i.e., $\curly{\a_j \colon j \in \supp(\x_l)}$ and $\curly{\a_j \colon j \in \supp(\x_{l'})}$, should be similar as well. 
\begin{definition}\label{def:similar}
    Let $0<\alpha,\beta<1$ be two predefined threshold parameters. We say that two atoms $\a_j$ and $\a_{j'}$ are \emph{similar} if
    \begin{align*}
        \corr(\a_j,\a_{j'}) = \frac{|\langle \a_j,\a_{j'} \rangle|}{\| \a_j \|_2 \| \a_{j'} \|_2} > \alpha.
    \end{align*} 
    Two measurements $\y_l$ and $\y_{l'}$ are \emph{similar} if
    \begin{align*}
        \corr(\y_l,\y_{l'}) = \frac{|\langle \y_l,\y_{l'} \rangle|}{\norm{\y_l}_2 \norm{\y_{l'}}_2} > \beta.
    \end{align*} 
\end{definition}
Let $\bar{\A}$ and $\bar{\Y}$ denote copies of $\A$ and $\Y$ with re-normalized columns, i.e., $\bar{\A}^T\bar{\A}$ and $\bar{\Y}^T\bar{\Y}$ contain the pairwise correlation of all atoms and measurement vectors. We define
\begin{align} \label{eq:Cbeta}
\C_{\beta}\in\{0,1\}^{L\times L}, && (C_{\beta})_{j,k}=\begin{cases}
1 & (\bar{Y}^T \bar{Y})_{j,k}>\beta\\
0 & (\bar{Y}^T \bar{Y})_{j,k}\leq\beta
\end{cases},
\end{align}
and $\C_{\alpha} \in \R^{N \times N}$ accordingly where $\bar{\Y}$ is replaced with $\bar{\A}$. Consequently, $\C_{\beta}$ is $1$ for each pair of similar measurement vectors and $0$ otherwise. A naive design for $\C$ is then
\begin{align} \label{eq:C}
\C_{\text{simple}} = \C_{\beta} \otimes (\mathbb{1}-\C_{\alpha}) \in \R^{NL \times NL}.
\end{align}
Now $C_{(j,l),(j',l')}$ is $1$ whenever the measurement vectors $\y_l$ and $\y_{l'}$ are similar but the chosen atoms $\a_{j}$ and $\a_{j'}$ are not.

Constructing $\C$ as in \eqref{eq:C}, \eqref{eq:pen2} penalizes whenever two measurement vectors $\y_l$ and $\y_{l'}$ are similar but the atoms of $\A$ activated by $\x_l$ and $\x_{l'}$ are not. It does not penalize activation of similar atoms $\a_{j}$ and $\a_{j'}$ by a single $\x_l$, i.e., in the same measurement vector $\y_l$. Numerical experiments show that blurring effects in the reconstruction are a consequence. To overcome the problem, we instead construct $\C$ as a composition of $L^2$ blocks of size $N\times N$ in the following way:

\begin{center}
	\begin{tikzpicture}[scale = 0.35]
	\node[left] at (0,1.05) {$\C=$};
	
	\draw[blue,thick] (0,-1.1) rectangle (1,-0.1);
	\draw[blue,thick] (1.1,-1.1) rectangle (2.1,-0.1);
	\draw[blue,thick] (2.2,-1.1) rectangle (3.2,-0.1);
	\draw[red,very thick,dotted] (3.3,-1.1) rectangle (4.3,-0.1);	
	
	\draw[blue,thick] (0,0) rectangle (1,1);
	\draw[blue,thick] (1.1,0) rectangle (2.1,1);
	\draw[red,very thick,dotted] (2.2,0) rectangle (3.2,1);
	\draw[blue,thick] (3.3,0) rectangle (4.3,1);
	
	\draw[blue,thick] (0,1.1) rectangle (1,2.1);
	\draw[red,very thick,dotted] (1.1,1.1) rectangle (2.1,2.1);
	\draw[blue,thick] (2.2,1.1) rectangle (3.2,2.1);
	\draw[blue,thick] (3.3,1.1) rectangle (4.3,2.1);
	
	\draw[red,very thick,dotted] (0,2.2) rectangle (1,3.2);
	\draw[blue,thick] (1.1,2.2) rectangle (2.1,3.2);
	\draw[blue,thick] (2.2,2.2) rectangle (3.2,3.2);
	\draw[blue,thick] (3.3,2.2) rectangle (4.3,3.2);
	
	\draw[green!50!black,very thick,dash dot] (0,3.2) -- (4.3,-1.1);
	
	\draw[blue,thick] (5.5,-1.1) rectangle (6.5,-0.1);
	\draw[blue,thick] (6.6,-1.1) rectangle (7.6,-0.1);
	\draw[blue,thick] (7.7,-1.1) rectangle (8.7,-0.1);
	\draw[red,very thick,dotted] (14.8,-1.1) rectangle (15.8,-0.1);	
	
	\draw[blue,thick] (5.5,0) rectangle (6.5,1);
	\draw[blue,thick] (6.6,0) rectangle (7.6,1);
	\draw[red,very thick,dotted] (13.7,0) rectangle (14.7,1);
	\draw[blue,thick] (8.8,0) rectangle (9.8,1);
	
	\draw[blue,thick] (5.5,1.1) rectangle (6.5,2.1);
	\draw[red,very thick,dotted] (12.6,1.1) rectangle (13.6,2.1);
	\draw[blue,thick] (7.7,1.1) rectangle (8.7,2.1);
	\draw[blue,thick] (8.8,1.1) rectangle (9.8,2.1);
	
	\draw[red,very thick,dotted] (11.5,2.2) rectangle (12.5,3.2);
	\draw[blue,thick] (6.6,2.2) rectangle (7.6,3.2);
	\draw[blue,thick] (7.7,2.2) rectangle (8.7,3.2);
	\draw[blue,thick] (8.8,2.2) rectangle (9.8,3.2);
	
	\draw[green!50!black,very thick,dash dot] (17.5,3.2) -- (21.8,-1.1);
	
	\node[right] at (10,1.05) {$+$};
	\node[right] at (16,1.05) {$+$};
	\node[right] at (4.2,1.05) {$=$};
	\node[right] at (4.3,-2.5) {$=\textcolor{blue}{\C_1}+\textcolor{red}{\C_2}+\textcolor{green!50!black}{\C_3}$.};
	\end{tikzpicture}
\end{center}
Here, $\C_1$ is the penalization for pairwise different measurement vectors, $\C_2$ compares each column of the solution with itself and $\C_3$ compares each element with itself. Since the comparison of each element with itself does not give any structural information, we set $\C_3=\0$. Pairwise different measurements can be treated as in the example above and thus
\begin{align} \label{eq:C1}
	\C_1=(\C_{\beta}-\id) \otimes (\mathbb{1}-\C_{\alpha}),
\end{align} 
where using $(\C_{\beta}-\id)$ instead of $\C_\beta$ sets $\C_1$ to zero on the diagonal blocks. For $\C_2$ we use the contrary heuristic: whenever two elements in the same column $\x_l$ of $\X$ are non-zero, the corresponding activated atoms in $\A$ should not be similar as otherwise one of them would be redundant. We set 
\begin{align} \label{eq:C2}
	\C_2=\id \otimes(\C_{\alpha}-\id),
\end{align} 
where subtracting the identity sets $\C_2$ to zero on the main diagonal. This construction of $\C$ exhibits sound performance in numerical experiments, see Section \ref{sec:Numerics}.

\begin{remark} \label{rem:Kronecker}
	The Kronecker form of the above matrices allows fast matrix-vector multiplication. This becomes important later on as $\C$ is of high dimension and multiplication by $\C$ is a fundamental operation for our numerical simulations.
\end{remark}

\noindent
\textbf{Choosing $\alpha$ and $\beta$:} 
When using the above construction method, we need to choose suitable parameters $\alpha$ and $\beta$ controlling the shape of $\C$. As in Proposition \ref{thm:Equivalence}, we assume here that each column of $\X$ contains at most one non-zero entry per active elementary structure, cf.\ Remark \ref{rem:Assumption}. 

\noindent
A meaningful parameter $\alpha$ can be directly deduced from prior knowledge on the concrete application. See Section \ref{sec:Numerics} for examples on how these priors can look like and how $\alpha$ is derived from them. If $\alpha$ is given, $\beta$ can be estimated from $\Y$ and $\A$. To this end, let $\y,\y'$ be two columns of $\Y$ and $\x,\x'$ the corresponding columns of $\X$, i.e., $\y = \A\x$ and $\y' = \A\x'$. Recall from \eqref{eq:C} that at this point we are interested in determining $\beta$ such that, by Definition \ref{def:similar}, $\y$ and $\y'$ are similar only if $\x$ and $\x'$ activate similar atoms. For simplicity, we assume here that $\A$ has unit norm columns; the argument can be easily adapted to the general case. Assume $s$ elementary structures are active in $\X$ where we do not need to know $s$. We introduce copies $\xxi,\xxi' \in \R^s$ of $\x,\x'$ which are reduced to the support of $\x,\x'$ and re-ordered such that $\xi_k$ and $\xi'_k$ belong to the same elementary structure, for all $k \in [s]$. If there is no noise on the measurements, we can write $\y$ and $\y'$ as
\begin{align*}
    \y=\sum\limits_{k=1}^s \xi_k \a_{j_k} \quad \text{ and } \quad
    \y'=\sum\limits_{k=1}^s \xi'_k \a_{j'_k},
\end{align*}
where $\a_{j_k},\a_{j'_k}$ are the atoms activated by the $k$-th structure in $\x$ resp.\ $\x'$. Hence, the inner product of $\y$ and $\y'$ is
\begin{align}
\label{eq:Decomposition}
    \langle \y,\y'\rangle &=
    \sum\limits_{k=1}^s \xi_k \xi'_k  \langle \a_{j_k}, \a_{j'_k} \rangle +
    \sum\limits_{\substack{k,k'=1\, \\ k\neq k'}}^s \xi_k \xi'_{k'} \langle \a_{j_k}, \a_{j'_{k'}}\rangle.
\end{align}
The first term represents the correlation of each of the $s$ structures with itself, while the second represents interference between different structures. Note that the first sum only contains inner products of atoms $\a_{j_k}$ and $\a_{j'_k}$ activated by the same structure. Since we consider here the situation that activated atoms are similar, by Definition \ref{def:similar}, the absolute value of these inner products is limited from below by $\alpha$. If we assume that 
\begin{enumerate}
    \item[\textbf{(A)}] there are no sudden phase shifts in the elementary structures of our ground-truth $\X$,\\
    i.e., $\sign(\xi_k) = \sign(\xi'_k)$, for all $k \in [s]$,
\end{enumerate}
then all terms in the first sum are positive. If, in addition\footnote{Both assumptions, (A) and (B), are mild and often hold in applications. Assumption (B), e.g., holds whenever the elementary structures are separated or destructive interference is occurring.},
\begin{enumerate}
    \item[\textbf{(B)}] the interference between elementary structures is bounded by $0 < \eps \ll \alpha \langle \xxi, \xxi' \rangle$,
\end{enumerate}
we know from \eqref{eq:Decomposition} that
\begin{align}\label{eq:beta1}
    \langle \y, \y' \rangle \ge \alpha \langle \xxi, \xxi' \rangle - \eps
    \quad \text{ and } \quad  
    \norm{\y}^2 \le \norm{\xxi}_2^2 + \varepsilon.
\end{align}
For the correlation of the two measurements we thus obtain
\begin{align}\label{eq:beta2}
\begin{split}
    \corr(\y,\y')&=\frac{\langle\y,\y'\rangle}{\|\y\|_2\|\y'\|_2}
    \geq\frac{\alpha \langle \xxi, \xxi' \rangle - \eps}{\sqrt{\|\xxi\|_2^{2}+\eps}\sqrt{\|\xxi'\|_2^{2}+\eps}}
    \\&\approx\alpha \, \corr(\xxi,\xxi'),
\end{split}
\end{align}
such that we can choose $\beta = \alpha \corr(\xxi,\xxi')$. Note that $\corr(\xxi,\xxi')$ is close to $1$ whenever the entries of $\X$ do hardly vary along elementary structures, and it is small for highly fluctuating entries along elementary structures. In applications the (expected) correlation $\corr(\xxi,\xxi')$ should be given approximately.
Let us mention that the above assumptions are not required to hold for arbitrary $\y,\y'$ and elementary structures. To have a reliable heuristic parameter choice, it suffices if they apply to the majority of the measurements.

In the presence of noise, we decrease the estimate \eqref{eq:beta2} relative to the noise level. By this the measurements are still considered similar even if they are corrupted by noise. The simulations in Section \ref{sec:Numerics} show efficacy of the heuristic.

Let us finally mention that while testing in numerical experiments the performance of \eqref{eq:LASSO_pos_SinglePenalty} with $\C$ as constructed above, it turned out that depending on the concrete problem setting the tuning of $\lambda$ is challenging. The main problem appears to be the different contributions of $\C_1$ and $\C_2$ in \eqref{eq:C1} and \eqref{eq:C2} to the regularizing function $\Rc_\C$. While $\C_1$ compares pairwise different measurements and enforces global structure, $\C_2$ compares each column of the solution with itself enforcing sparsity along the columns. Consequently, it is beneficial to additionally balance between $\C_1$ and $\C_2$ by introducing parameters $\lambda_1,\lambda_2 > 0$ and defining $\C_\LAmbda = \lambda_1 \C_1 + \lambda_2 \C_2$, for $\LAmbda = (\lambda_1,\lambda_2)$. Then, the program \eqref{eq:LASSO_pos_SinglePenalty} becomes
\begin{align} \label{eq:LASSO_pos}
    \min_{\Z \in \R_+^{N \times L}} \| \A\Z - \Y \|_F^2 + \lambda_1 \Rc_{\C_1} (\Z)+\lambda_2 \Rc_{\C_2} (\Z).
\end{align}
The preceding results --- Lemma \ref{lem:Lasso}, Propostion \ref{thm:Geometry}, and Corollary \ref{cor:Geometry} --- extend in a straight-forward way to \eqref{eq:LASSO_pos} by using linearity of $\Rc$ in $\C$, i.e., $\lambda_1 \Rc_{\C_1} (\Z)+\lambda_2 \Rc_{\C_2} (\Z) = \Rc_{\C_\LAmbda}(\Z)$.

\begin{remark} \label{rmk:parameterchoice}
    Let us briefly comment on why $\lambda_1$ and $\lambda_2$ may notably differ in applications. Recall the proof of Proposition \ref{thm:Equivalence}. In particular, the bounds in inequality \eqref{eq:proofInequality} are dominated by the influence of $\C_1$. Considering both matrix parts separately, one could refine the bound to be
    \begin{align*}
        1 \le \vv{\X_k}^T \C_1 \vv{\X_j} \le L^2 \quad \text{and} \quad
        0 \le \vv{\X_k}^T \C_2 \vv{\X_j} \le L.
    \end{align*}
    The second term penalizes structures that are overlapping or close to each other. Since this hardly happens throughout all measurements, we expect the upper bound to be overpessimistic. In applications, one would rather have $\mathcal{O}(1)$ than $\mathcal{O}(L)$. Considering $\C_1$, i.e., the first inequality, however, a scaling of $L^2$ seems realistic whenever the structures are active throughout all measurements. Hence, the influence of $\C_1$ on the penalty term is up to $L^2$-times larger compared to $\C_2$ which can be compensated by scaling $\lambda_2$ accordingly. In addition, $\lambda_1$ balances the structural sparsity and must be chosen smaller the more structures appear in the data while $\lambda_2$ handles the sparsity within single structures and is mostly independent of the number of active elementary structures. As Section \ref{sec:Numerics} shows, a parameter setup with very small $\lambda_1 = \mathcal{O}(L^{-2})$ and $\lambda_2 \gg 1$ is not unusual.
\end{remark}

	
\bibliographystyle{IEEEtranS}
\bibliography{IEEEabrv,references}

\begin{thebibliography}{10}
\providecommand{\url}[1]{#1}
\csname url@samestyle\endcsname
\providecommand{\newblock}{\relax}
\providecommand{\bibinfo}[2]{#2}
\providecommand{\BIBentrySTDinterwordspacing}{\spaceskip=0pt\relax}
\providecommand{\BIBentryALTinterwordstretchfactor}{4}
\providecommand{\BIBentryALTinterwordspacing}{\spaceskip=\fontdimen2\font plus
\BIBentryALTinterwordstretchfactor\fontdimen3\font minus
  \fontdimen4\font\relax}
\providecommand{\BIBforeignlanguage}[2]{{%
\expandafter\ifx\csname l@#1\endcsname\relax
\typeout{** WARNING: IEEEtranS.bst: No hyphenation pattern has been}%
\typeout{** loaded for the language `#1'. Using the pattern for}%
\typeout{** the default language instead.}%
\else
\language=\csname l@#1\endcsname
\fi
#2}}
\providecommand{\BIBdecl}{\relax}
\BIBdecl

\bibitem{aharon2006k}
M.~Aharon, M.~Elad, and A.~Bruckstein, ``K-svd: An algorithm for designing
  overcomplete dictionaries for sparse representation,'' \emph{IEEE
  Transactions on signal processing}, vol.~54, no.~11, pp. 4311--4322, 2006.

\bibitem{Angelosante09}
D.~Angelosante, G.~B. Giannakis, and E.~Grossi, ``Compressed sensing of
  time-varying signals,'' in \emph{2009 16th International Conference on
  Digital Signal Processing}, 2009, pp. 1--8.

\bibitem{baraniuk2010model}
R.~G. Baraniuk, V.~Cevher, M.~F. Duarte, and C.~Hegde, ``Model-based
  compressive sensing,'' \emph{{IEEE} Trans. Inf. Theory}, vol.~56, no.~4, pp.
  1982--2001, 2010.

\bibitem{Baron06}
D.~Baron, M.~B. Wakin, M.~F. Duarte, S.~Sarvotham, and R.~G. Baraniuk,
  ``Distributed compressed sensing,'' \emph{IEEE Trans. Information Theory},
  vol.~52, no.~12, pp. 5406--5425, 2006.

\bibitem{Bossmann12}
F.~Boßmann, G.~Plonka, T.~Peter, O.~Nemitz, and T.~Schmitte, ``Sparse
  deconvolution methods for ultrasonic ndt,'' \emph{Journal of Nondestructive
  Evaluation}, vol.~31, no.~3, pp. 225--244, 2012.

\bibitem{Bossmann15}
F.~Bossmann and J.~Ma, ``Asymmetric chirplet transform for sparse
  representation of seismic data,'' \emph{Geophysics}, vol.~80, no.~6, pp.
  WD89--WD100, 2015.

\bibitem{Bossmann16}
------, ``Asymmetric chirplet transform—part 2: Phase, frequency, and chirp
  rate,'' \emph{Geophysics}, vol.~81, no.~6, pp. V425--V439, 2016.

\bibitem{candes2006robust}
E.~J. Cand{\`e}s, J.~Romberg, and T.~Tao, ``Robust uncertainty principles:
  Exact signal reconstruction from highly incomplete frequency information,''
  \emph{IEEE Trans. Inf. Theory}, vol.~52, no.~2, pp. 489--509, 2006.

\bibitem{candes2006near}
E.~J. Cand\`es and T.~Tao, ``Near-optimal signal recovery from random
  projections: Universal encoding strategies?'' \emph{IEEE Trans. Inf. Theory},
  vol.~52, no.~12, pp. 5406--5425, 2006.

\bibitem{chambolle2011first}
A.~Chambolle and T.~Pock, ``A first-order primal-dual algorithm for convex
  problems with applications to imaging,'' \emph{Journal of mathematical
  imaging and vision}, vol.~40, no.~1, pp. 120--145, 2011.

\bibitem{chandrasekaran2010convex}
V.~Chandrasekaran, B.~Recht, P.~A. Parrilo, and A.~S. Willsky, ``The convex
  algebraic geometry of linear inverse problems,'' in \emph{2010 48th Annual
  Allerton Conference on Communication, Control, and Computing
  (Allerton)}.\hskip 1em plus 0.5em minus 0.4em\relax IEEE, 2010, pp. 699--703.

\bibitem{Chen06}
J.~Chen and X.~Huo, ``Theoretical results on sparse representations of
  multiple-measurement vectors,'' \emph{IEEE Trans. Signal Processing},
  vol.~54, no.~12, pp. 4634--4643, 2006.

\bibitem{Cotter05}
S.~F. Cotter, B.~D. Rao, K.~Engan, and K.~Kreutz{-}Delgado, ``Sparse solutions
  to linear inverse problems with multiple measurement vectors,'' \emph{IEEE
  Trans. Signal Process.}, vol.~53, no.~7, pp. 2477--2488, 2005.

\bibitem{daubechies2004iterative}
I.~Daubechies, M.~Defrise, and C.~De~Mol, ``An iterative thresholding algorithm
  for linear inverse problems with a sparsity constraint,''
  \emph{Communications on Pure and Applied Mathematics}, vol.~57, no.~11, pp.
  1413--1457, 2004.

\bibitem{davies12}
M.~Davies and Y.~Eldar, ``Rank awareness in joint sparse recovery.'' \emph{IEEE
  Trans. Information Theory}, vol.~58, no.~2, pp. 1135--1146, 2012.

\bibitem{donoho2006compressed}
D.~L. Donoho, ``Compressed sensing,'' \emph{IEEE Trans. Inf. Theory}, vol.~52,
  no.~4, pp. 1289--1306, 2006.

\bibitem{foucart2013mathematical}
S.~Foucart and H.~Rauhut, \emph{A Mathematical Introduction to Compressive
  Sensing}.\hskip 1em plus 0.5em minus 0.4em\relax Birkh\"{a}user Basel, 2013.

\bibitem{Heins2015}
P.~Heins, M.~Moeller, and M.~Burger, ``Locally sparse reconstruction using the
  {$\ell^{1,\infty}$}-norm.'' \emph{Inverse Problems and Imaging}, vol.~9,
  no.~4, pp. 1093--1137, 2015.

\bibitem{Huang2011}
J.~Huang, T.~Zhang, and D.~Metaxas, ``Learning with structured sparsity.''
  \emph{Journal of Machine Learning Research}, vol.~12, no.~11, pp. 3371--3412,
  2011.

\bibitem{huang2010benefit}
J.~Huang and T.~Zhang, ``The benefit of group sparsity,'' \emph{The Annals of
  Statistics}, vol.~38, no.~4, pp. 1978--2004, 2010.

\bibitem{Katz87}
S.~M. Katz, ``Estimation of probabilities from sparse data for the language
  model component of a speech recognizer,'' \emph{{IEEE} Trans. Acoust. Speech
  Signal Process.}, vol.~35, no.~3, pp. 400--401, 1987.

\bibitem{Lee2014}
S.~Lee, Y.~Liao, M.~Seo, and Y.~Shin, ``Structural change in sparsity.''
  \emph{arXiv preprint arXiv:1411.3062}, 2014.

\bibitem{Lian2019}
L.~Lian, A.~Liu, and V.~Lau, ``Exploiting dynamic sparsity for downlink
  fdd-massive mimo channel tracking.'' \emph{IEEE Trans. Signal Processing},
  vol.~67, no.~8, pp. 2007--2021, 2019.

\bibitem{Lu12}
X.~Lu, H.~Yuan, P.~Yan, Y.~Yuan, and X.~Li, ``Geometry constrained sparse
  coding for single image super-resolution,'' \emph{2012 IEEE Conf. Computer
  Vision and Pattern Recognition}, pp. 1648--1655, 2012.

\bibitem{Lustig07}
M.~Lustig, D.~Donoho, and J.~M. Pauly, ``Sparse mri: The application of
  compressed sensing for rapid mr imaging,'' \emph{Magnetic Resonance in
  Medicine}, vol.~58, no.~6, pp. 1182--1195, 2007.

\bibitem{peleg2012exploiting}
T.~Peleg, Y.~C. Eldar, and M.~Elad, ``Exploiting statistical dependencies in
  sparse representations for signal recovery,'' \emph{{IEEE} Trans. Signal
  Process.}, vol.~60, no.~5, pp. 2286--2303, 2012.

\bibitem{Rao98}
B.~D. Rao, ``Signal processing with the sparseness constraint,'' in
  \emph{Proceedings of the 1998 {IEEE} International Conference on Acoustics,
  Speech and Signal Processing, {ICASSP} '98, Seattle, Washington, USA, May
  12-15, 1998}.\hskip 1em plus 0.5em minus 0.4em\relax {IEEE}, 1998, pp.
  1861--1864.

\bibitem{rao:2014}
X.~Rao and V.~K.~N. Lau, ``Distributed compressive {CSIT} estimation and
  feedback for {FDD} multi-user massive {MIMO} systems,'' \emph{{IEEE} Trans.
  Signal Process.}, vol.~62, no.~12, pp. 3261 -- 3271, 2014.

\bibitem{rick2017one}
J.~Rick~Chang, C.-L. Li, B.~Poczos, B.~Vijaya~Kumar, and A.~C.
  Sankaranarayanan, ``One network to solve them all--solving linear inverse
  problems using deep projection models,'' in \emph{Proceedings of the IEEE
  International Conference on Computer Vision}, 2017, pp. 5888--5897.

\bibitem{Sidky08}
E.~Y. Sidky and X.~Pan, ``Image reconstruction in circular cone-beam computed
  tomography by constrained, total-variation minimization,'' \emph{Physics in
  medicine and biology}, vol.~53, no.~17, pp. 4777--4807, 2008.

\bibitem{Starck10}
J.-L. Starck, F.~Murtagh, and J.~Fadili, \emph{Sparse Image and Signal
  Processing: Wavelets, Curvelets, Morphological Diversity}.\hskip 1em plus
  0.5em minus 0.4em\relax Cambridge University Press, 2010.

\bibitem{Tropp06b}
J.~A. Tropp, ``Algorithms for simultaneous sparse approximation: part ii:
  Convex relaxation,'' \emph{Signal Process.}, vol.~86, no.~3, pp. 589 -- 602,
  2006.

\bibitem{Tropp06a}
J.~A. Tropp, A.~C. Gilbert, and M.~J. Strauss, ``Algorithms for simultaneous
  sparse approximation. part {I:} greedy pursuit,'' \emph{Signal Process.},
  vol.~86, no.~3, pp. 572--588, 2006.

\bibitem{Vaswani2010}
N.~Vaswani and W.~Lu, ``Modified-cs: Modifying compressive sensing for problems
  with partially known support.'' \emph{IEEE Trans. Signal Processing},
  vol.~58, no.~9, pp. 4595--4607, 2010.

\bibitem{Wipf07}
D.~Wipf and B.~D. Rao, ``An empirical bayesian strategy for solving the
  simultaneous sparse approximation problem,'' \emph{IEEE Trans. Signal
  Process.}, vol.~55, no.~7, pp. 3704--3716, 2007.

\bibitem{Wright09}
J.~Wright, A.~Y. Yang, A.~Ganesh, S.~S. Sastry, and Y.~Ma, ``Robust face
  recognition via sparse representation,'' \emph{IEEE Trans. Pattern Anal.
  Mach. Intell.}, vol.~31, no.~2, pp. 210--227, 2009.

\bibitem{wu:2014}
Y.~Wu, Y.-J. Zhu, Q.-Y. Tang, C.~Zou, W.~Liu, R.-B. Dai, X.~Liu, E.~X. Wu,
  L.~Ying, and D.~Liang, ``Accelerated {MR} diffusion tensor imaging using
  distributed compressed sensing,'' \emph{Magnetic Resonance in Medicine},
  vol.~71, no.~2, pp. 764 -- 772, 2014.

\bibitem{yu2012bayesian}
L.~Yu, H.~Sun, J.-P. Barbot, and G.~Zheng, ``Bayesian compressive sensing for
  cluster structured sparse signals,'' \emph{Signal Process.}, vol.~92, no.~1,
  pp. 259--269, 2012.

\bibitem{yuan2006model}
M.~Yuan and Y.~Lin, ``Model selection and estimation in regression with grouped
  variables,'' \emph{Journal of the Royal Statistical Society: Series B
  (Statistical Methodology)}, vol.~68, no.~1, pp. 49--67, 2006.

\bibitem{sblWebpage}
\BIBentryALTinterwordspacing
Z.~Zhang, ``(5) {SBL} ({S}parse {B}ayesian {L}earning),'' ver. 1.1
  (02/12/2011), (accessed 09/02/2021). [Online]. Available:
  \url{http://dsp.ucsd.edu/~zhilin/Software.html}
\BIBentrySTDinterwordspacing

\bibitem{Zhang10}
Z.~Zhang and B.~D. Rao, ``Sparse signal recovery in the presence of correlated
  multiple measurement vectors,'' in \emph{Proceedings of the {IEEE}
  International Conference on Acoustics, Speech, and Signal Processing,
  {ICASSP} 2010, 14-19 March 2010, Sheraton Dallas Hotel, Dallas, Texas,
  {USA}}, 2010, pp. 3966--3989.

\bibitem{zhang2013extension}
------, ``Extension of sbl algorithms for the recovery of block sparse signals
  with intra-block correlation,'' \emph{IEEE Trans. Signal Process.}, vol.~61,
  no.~8, pp. 2009--2015, 2013.

\bibitem{Zheng11}
H.~Zheng, J.~Bu, C.~Chen, C.~Wang, L.~Zhang, G.~Qiu, and D.~Cai, ``Graph
  regularized sparse coding for image representation,'' \emph{IEEE Trans. Image
  Processing}, vol.~20, no.~5, pp. 1327--1336, 2011.

\bibitem{zhu2011sparsity}
H.~Zhu, G.~Leus, and G.~B. Giannakis, ``Sparsity-cognizant total least-squares
  for perturbed compressive sampling,'' \emph{IEEE Transactions on Signal
  Processing}, vol.~59, no.~5, pp. 2002--2016, 2011.

\bibitem{Ziniel13a}
J.~Ziniel and P.~Schniter, ``Dynamic compressive sensing of time-varying
  signals via approximate message passing,'' \emph{IEEE Trans. Signal
  Process.}, vol.~61, no.~21, pp. 5270--5284, 2013.

\bibitem{Ziniel13}
------, ``Efficient high-dimensional inference in the multiple measurement
  vector problem,'' \emph{{IEEE} Trans. Signal Process.}, vol.~61, no.~2, pp.
  340--354, 2013.

\end{thebibliography}
	
\end{document}